\documentclass{article}
\usepackage{mathpazo}
\usepackage{graphicx} 
\usepackage{amsmath, amssymb, amsthm, amsfonts}
\usepackage{comment}
\usepackage[margin=1in]{geometry}
\usepackage{color}
\usepackage{hyperref}
\usepackage{float}

\usepackage{thmtools}
\declaretheoremstyle[headfont=\kpfonts]{normalhead}

\newtheorem{theorem}{Theorem}[section]
\newtheorem{proposition}{Proposition}[section]
\newtheorem{lemma}{Lemma}[section]

\theoremstyle{definition}

\newtheorem{remark}{Remark}[section]

\renewcommand{\Re}{\operatorname{Re}}
\renewcommand{\Im}{\operatorname{Im}}

\title{Fokas method for linear convection-diffusion equation with time-dependent coefficients and its extension to other evolution equations}

\author{%
\begin{minipage}{0.45\textwidth}
  \centering
  Konstantinos Kalimeris \\
  \small Mathematics Research Center, Academy of Athens, Athens, Greece
\end{minipage}
\hfill
\begin{minipage}{0.45\textwidth}
  \centering
  Türker Özsarı \\
  \small Department of Mathematics, Bilkent University, Ankara, Turkey
\end{minipage}
}

\date{}

\begin{document}

\maketitle
\begin{abstract}
In this paper, we study a linear convection-diffusion equation with time-dependent coefficients on a bounded interval, motivated by associated physical problems. 
We apply and adapt the Unified Transform Method (UTM), a.k.a. Fokas Method, which handles both time-varying coefficients and nonzero boundary data, to obtain an explicit integral formula for the solution.
Next, we study well-posedness of the model in fractional Sobolev spaces and prove spatial and temporal regularity estimates, showing that the smoothing effect of the heat operator is still prevalent even when coefficients depend on time.  Finally, we extend this approach to obtain the solution for several evolution equations with time-dependent coefficients, in one space variable.
\end{abstract}



\tableofcontents
\section{Introduction}

We study the well-posedness for the convection-diffusion (a.k.a advection-diffusion) equation with time dependent coefficients on a finite interval $\Omega= (0,L)\ni x$ subject to inhomogeneous boundary manipulations:
\begin{equation}\label{maineq}
	\begin{aligned}
		u_t = A(t)u; \quad t\in (0,T);\quad 
		u(0) = u_0;\quad
		\gamma_0 u=(g_0,h_0),
	\end{aligned}
\end{equation}
where $A(t)$ is the differential operator formally given by \begin{equation}\label{A_t}A(t)u=b(t)u_{xx}+c(t)u_x,\end{equation}
where $b(t)\ge b_0$ for some $b_0>0$ and $c(t)\in\mathbb{R}$ for all $t\ge 0$, and $\gamma_0$ is the (Dirichlet) boundary operator formally defined by $\gamma_0 u = (u|_{x=0},u|_{x=L}).$ This is achieved by first deriving an integral representation of the solution and obtaining the associated sharp regularity estimates.

\maketitle
Apart from the mathematical curiosity,  convection-diffusion model on a finite region with time dependent coefficients and inhomogeneous boundary conditions,  is a core model for real-world systems undergoing transport, diffusion, and control, with wide applications. The coefficient $b(t)$ helps to describe random spreading or smoothing and it allows diffusivity to alternate due to physical changes.  On the other hand, the convection coefficient $c(t)$ captures non-steady or fluctuating transport processes. This is essential for simulating real-world systems where the velocity of the medium is not constant but evolves dynamically due to external factors or control mechanisms.  Thus, incorporating a time-dependent drift coefficient transforms the model from a simplistic approximation into a flexible, dynamic tool capable of representing the complex, time-varying nature of transport in physical, environmental, and engineered systems.  Inhomogeneous boundary conditions reflect real-world inputs (e.g., pollutants injected at endpoints), sinks, and interfaces, and bring the model closer to applications incorporating the temperature or concentration changes prescribed at boundary (reactor walls or riverbanks are just some examples).  The behavior of the system is dominated by the interaction between internal dynamics and boundary manipulations, leading to rich transient behavior and nontrivial steady states.  Finally, finiteness of the domain allows us to simulate physical mediums such as pipes, rods, and chambers.  

Certain forms of diffusion equations with time dependent coefficients have been solved semi-analytically using the idea of the generalized integral transform method and/or the method of heat potentials, see \cite{Itkin22}, \cite{ItLiDm},\cite{Car21} and the references therein. These problems, although similar, are structurally distinct from the ones we consider in the current manuscript, and the related methodology does not seem efficient for establishing sharp regularity results for \eqref{maineq}.

Problems associated to \eqref{maineq} especially attracted the attention of the inverse problems community, see for instance \cite{Cannon63} (also \cite{Cannon91}) where the author focuses on the determination of the conductivity of a medium if the conductivity was known a priori to be a
function of time only under the assumption that the Dirichlet boundary inputs are non-zero temporal functions.  See also \cite{Ivanchov1994547} where the author establishes existence and uniqueness conditions for the inverse problem of simultaneously determining time-dependent thermal conductivity and specific heat capacity in a one-dimensional heat equation, using analytic transformations and fixed-point arguments to prove solvability under suitable smoothness and boundary assumptions. Similarly, in \cite{Ivanchov1994547}, the author analyzes the inverse problem of simultaneously determining two time-dependent coefficients (the leading and lower-order terms) in a one-dimensional parabolic equation, proving local existence via Schauder’s fixed-point theorem and establishing global uniqueness under suitable boundary and regularity conditions.  There are also some studies for higher order diffusion equations, see for instance \cite{CLI21} where authors are reconstructing the time-dependent thermal grooving coefficient.

There has also been numerical efforts in this direction. For instance, \cite{Dehghan05} introduces and compares several finite difference schemes to identify a time-dependent coefficient in a one-dimensional parabolic PDE from overspecified boundary data, showing that the approaches can stably and accurately reconstruct both the solution and the coefficient from an extra measurement. \cite{HUNTUL2017} develops a numerical optimization method to recover time-dependent thermal conductivity in the one-dimensional heat equation from boundary heat flux data, reformulating the inverse problem as a nonlinear least-squares minimization, and demonstrating accuracy and stability through test examples. Similar concerns for the inverse problems associated with the convection-diffusion equation have also been addressed, see for instance \cite{Fister06} for a variational approach to the problem of 	
determination of a sorption parameter from soil column experiments.  These are only small fraction of the work on this subject, there is a vast literature both from analysis and numerical perspectives on inverse problems associated with diffusion equations with time dependent coefficients.  

Despite its relevance, this model poses analytical challenges due to its non-autonomous structure and boundary driven dynamics. Classical techniques (separation of variables, eigenfunction expansions) are limited in non-autonomous PDE settings. Traditional transform methods (Laplace transform in time) usually assume constant coefficients, and not the best approach for higher order PDEs, especially involving multiple spatial derivatives.  

In this paper, we apply and extend the Unified Transform Method (UTM), originally developed for constant-coefficient linear PDEs  \cite{FokasBook, DTV14} to analyze the initial-boundary value problem for a convection-diffusion equation on a finite interval with time dependent coefficients. The novelty lies in handling time-dependent diffusion and drift coefficients; inhomogeneous Dirichlet boundary data, and the derivation of explicit solution representations that eliminate unknown boundary traces using spectral symmetry and complex analysis.  Furthermore, we establish well-posedness results for the associated Cauchy and boundary value problems in fractional Sobolev spaces, obtaining both spatial and temporal regularity estimates which demonstrates the smoothing effect of the evolving heat semigroup, even in the presence of time-dependent coefficients under suitable sign and growth assumptions on them.  Although, we study only the direct problem in this manuscript, the inverse problem can also be handled, at least numerically, by using the Fokas method based control technique which we originally introduced for boundary control problems in \cite{KalOzs24}.  This is postponed to a future work; to our knowledge related work has been recently performed analytically and numerically, by employing the Fokas method for the case of constant coefficients, \cite{KalMinPal25, KalMin-pre}. As a side note, we only consider time dependent coefficients here.  For the case of the direct problem for spatially varying coefficients, see the recent work of \cite{Defa25}.  It should be noted that the case of spatially varying coefficients is not amenable to standard spatial Fourier transform and obtaining a solution representation is more challenging.

The unified transform method based well-posedness analysis of solutions of evolution PDEs in fractional Sobolev spaces was introduced a decade ago in a series of papers, see e.g., \cite{HM95},   \cite{FHM16}, and \cite{FHM17} for the earlier work in this direction. Later many other papers have been written by the same authors and others (see e.g., \cite{H23}, \cite{M23}, \cite{ATM24}, \cite{OY19}, and the references therein). Nowadays, the approach of UTM based well-posedness analysis of solutions of ibvps in special geometries such as the half-space or the finite line, is now accepted as one of the closest analogues of Fourier transform based analysis of solutions of initial value problems (IVPs). 

UTM was also used to study the well-posedness of the initial-boundary value problem (IBVP) associated with the diffusion equation with constant coefficients in fractional Sobolev spaces in \cite{HMY19}. The present paper extends the elegant theory presented in \cite{HMY19} in several directions: (i) coefficients are allowed to vary with time variable; (ii) the main equation involves multiple spatial derivatives; (iii) certain regularity estimates are boosted utilizing the parabolic smoothing effect of the heat operator. 

This paper is structured as follows: in Section \ref{deriv}, we establish the representation formula which will be used to define weak solutions of \eqref{maineq}. The same formula is also used for constructing numerical solutions in Section \ref{simul}. Section \ref{HLsetting} extends the representation formula to the half-line setting. In Section \ref{wpsection}, we decompose the linear PDE analysis into three simpler problems, two of which have to do with Cauchy problems which will later be associated with spatial extensions of initial and interior data (see Section \ref{Seccauchy}), and the third one is a reduced IBVP which will later be associated with modified boundary conditions obtained by subtracting the traces of solutions of aforementioned Cauchy problems at the boundary points (see Section \ref{Secredibvp}).  In the case of zero drift, local solutions are also global, however in the case of nonzero drift, this is no longer trivial and extra effort has to be given. This is done in Section \ref{Secglobal}.  In Section \ref{Secevol}, we extend the UTM formulation at the representational level to other evolution PDEs with time dependent coefficients. In Section \ref{Secfinal}, we make connection to the associated halfline problem regarding the well-posedness.  In addition, some comments are provided for the nonlinear problems within the same section. 

\section{Integral representation formula for weak solutions}\label{sec-sol-adv-dif}
\subsection{Derivation of representation formula}\label{deriv}
We first introduce the finite-line Fourier transform (FLFT) and its inverse as the pair of formulas below:

\begin{equation}
\hat{\varphi}(k):=\int_0^Le^{-ikx}\varphi(x)dx, \,k\in \mathbb{C};\quad \varphi(x):=\frac{1}{2\pi}\int_{-\infty}^\infty e^{ikx}\hat{\varphi}(k)dk, \,x\in (0,L).
\end{equation}

Note that spectral variable is allowed to be complex in the forward formula.  We apply FLFT to the main equation in \eqref{maineq}, which under integration by parts yields:

\begin{equation}
\begin{aligned}
	\hat{u}&_t(k,t) 
	=  b(t)\int_0^L e^{-ikx} u_{xx}(x,t)\,dx 
	+ c(t)\int_0^L e^{-ikx} u_{x}(x,t)\,dx \\
	&=  b(t)\left[h_1(t)e^{-ikL} - g_1(t) 
	+ ik h_0(t)e^{-ikL} - ik g_0(t)\right] \\
	&\quad + c(t)\left[h_0(t)e^{-ikL} - g_0(t)\right] 
	+ (-b(t)k^2 + ik c(t))\hat{u}(k,t)\\
	&=:\text{traces}(k,t)  + \big(-b(t)k^2 + ik c(t)\big)\hat{u}(k,t),
\end{aligned}
\end{equation} where $h_1 (t):= u_x(L,t)$ and $g_1(t):= u_x(0,t)$ are the unknown traces which have not been prescribed as data.
Rearranging the terms, we get
\begin{equation}\label{du_hat_dt}
\frac{d}{dt}\left[\hat{u}(k,t)e^{\omega(k,t)}\right] = e^{\omega(k,t)}\cdot \text{traces}(k,t),
\end{equation}where $\omega(k,t)=k^2B(t)-ikC(t)$ in which $B(t)$ and $C(t)$ are respectively anti-derivatives of $b(t)$ and $c(t)$ given by $B(t)=\int_0^tb(s)ds,\quad C(t)=\int_0^tc(s)ds, \quad t\ge 0.$  Integrating \eqref{du_hat_dt} over $[0,t]$, we get 
\begin{equation}\label{hat_u_first_form}
\hat{u}(k,t)=e^{-\omega(k,t)} \hat{u}_0(k) + e^{-\omega(k,t)}\int_0^te^{\omega(k,s)}\cdot \text{traces}(k,s)ds.
\end{equation}
For a given function $\psi:(0,T)\rightarrow \mathbb{R}$ of sufficient regularity, we define its $t$-transform with respect to a fixed map $f$ through the formula: $    \psi^f(k,t) := \int_0^tf(s)e^{\omega(k,s)}\psi(s)ds,\quad t\in (0,T).$
Using this notation, we can rewrite the last term in \eqref{hat_u_first_form} more explicitly as
\begin{equation}\label{traces}
\begin{aligned}
	\int_0^t&e^{\omega(k,s)}\cdot \text{traces}(k,s)ds \\
	= \,&e^{-ikL}h_1^b(k,t)-g_1^b(k,t)+ike^{-ikL}h_0^b(k,t)-ikg_0^b(k,t)+e^{-ikL}h_0^c(k,t)-g_0^c(k,t).
\end{aligned}
\end{equation} so that \eqref{hat_u_first_form} takes the following form, to which we will refer to as a global relation.

\begin{equation}\label{global_rel_1}
\begin{aligned}
	\hat{u}(k,t) \,=\, e^{-\omega(k,t)} \bigg[&\hat{u}_0(k) + e^{-ikL}h_1^b(k,t)-g_1^b(k,t)+ike^{-ikL}h_0^b(k,t)-ikg_0^b(k,t)\\
	&+e^{-ikL}h_0^c(k,t)-g_0^c(k,t)\bigg], \quad k\in\mathbb{C}.
\end{aligned}
\end{equation}

By applying the inverse FLFT to the global relation \eqref{global_rel_1} we obtain an integral representation of the solution which includes transforms of both known and unknown boundary traces:

\begin{align}
u(x,t)&=\frac{1}{2\pi} \int_\mathbb{R} e^{ikx-\omega(k,t)} \hat{u}_0(k) dk \notag \\
&- \frac{1}{2\pi} \int_\mathbb{R} e^{ikx-\omega(k,t)} \left[ ikg_0^b(k,t)+g_1^b(k,t)+g_0^c(k,t)    \right] dk \label{firstinverse}\\
&+\frac{1}{2\pi} \int_\mathbb{R} e^{-ik(L-x)-\omega(k,t)} \left[ ikh_0^b(k,t)+h_1^b(k,t)+h_0^c(k,t)    \right] dk.\notag
\end{align}

\begin{figure}
\centering
\includegraphics[width=0.45\linewidth]{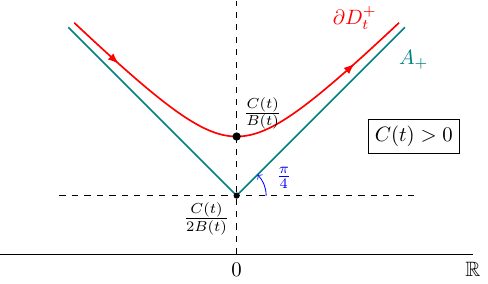} \qquad \includegraphics[width=0.45\linewidth]{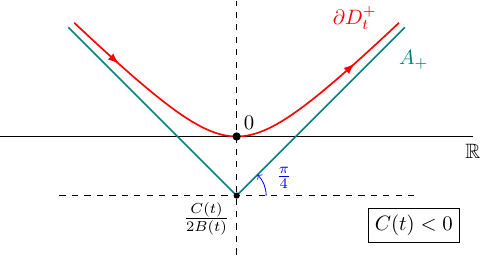}
\caption{The contours of integration $\partial D_t^+$.}
\label{fig:domain}
\end{figure}

To eliminate the unknowns we first introduce  a particular region in the complex plane via the assignment $D_t:=\{k\in\mathbb{C}\,|\,\text{Re}(\omega(k,t))<0\}$ and $\partial D_t$ as the boundary of $D_t$ in the complex plane. It's orientation is fixed so that as it is traversed $D_t$ remains at the left of it, see figure \ref{fig:domain}. The above definition implies that for every fixed $t\in(0,T)$ the contours $\partial D_t^\pm$ are suitable for following the standard procedure of the Fokas method. Namely, the second integral in \eqref{firstinverse} can be deformed to $\partial D_t^+$ and the third integral to $-\partial D_t^-$. Furthermore by defining 
\begin{equation}\label{def:HG}
G_1(k,t):= g_1^b(k,t)+g_0^c(k,t) ; \ \ \ H_1(k,t):= h_1^b(k,t)+h_0^c(k,t),
\end{equation}we obtain
\begin{align}\label{int_rep_1a}
u(x,t)&=\frac{1}{2\pi} \int_\mathbb{R} e^{ikx-\omega(k,t)} \hat{u}_0(k) dk \notag \\
&- \frac{1}{2\pi} \int_{\partial D_t^+} e^{ikx-\omega(k,t)} \left[ ikg_0^b(k,t)+G_1(k,t)    \right] dk \\
&-\frac{1}{2\pi} \int_{\partial D_t^-} e^{-ik(L-x)-\omega(k,t)} \left[ ikh_0^b(k,t)+H_1(k,t)    \right]  dk.\notag
\end{align}

An important step to eliminate unknowns $h_1^b$ and $g_1^b$ from the formula \eqref{hat_u_first_form}, we search for a map $k\mapsto \nu(k)$ that keeps the spectral input $\omega(k,t)$ invariant. To this end, we solve the equation $\omega(k,t)=\omega(\nu,t)$, which leads to two analytic solutions, one of which is the identity map, and the other is given by $\nu(k)=-k+i\frac{C(t)}{B(t)}, \quad k\in \mathbb{C}$; {where, for now, $C/B$ is assumed constant; this assumption is removed later in the well-posedness theory.} 

For convenience in computations we rewrite the global relation, by using \eqref{def:HG}, in the form 
\begin{equation}\label{global_rel_1a}
\begin{aligned}
	\hat{u}(k,t)e^{\omega(k,t)} &\,=\,\hat{u}_0(k) + e^{-ikL}H_1(k,t)-G_1(k,t)+ike^{-ikL}h_0^b(k,t)-ikg_0^b(k,t), \, k\in\mathbb{C}.
\end{aligned}
\end{equation}

Replacing $k$ with $\nu(k)$ leaves $H_1(k,t), \ G_1(k,t),\  h_0^b(k,t), \ g_0^b(k,t)$ invariant. Hence, under this mapping, the global relation \eqref{global_rel_1a}, yields
\begin{equation}\label{global_rel_2a}
\begin{aligned}
	\hat{u}&(\nu(k),t)e^{\omega(k,t)} \,=\,\hat{u}_0\left(\nu(k)\right) + e^{ikL+L\frac{C(t)}{B(t)}}H_1(k,t)-G_1(k,t)\\
	&-\left(ik+\frac{C(t)}{B(t)}\right)e^{ikL+L\frac{C(t)}{B(t)}}h_0^b(k,t) +\left(ik+\frac{C(t)}{B(t)}\right)g_0^b(k,t) , \, k\in\mathbb{C}.
\end{aligned}
\end{equation}

In what follows, we solve \eqref{global_rel_1a} and \eqref{global_rel_2a} for $H_1$ and $G_1$, and substitute the result into the integral representation \eqref{int_rep_1a} to obtain:

\begin{equation}\label{sol-rep-1}
\begin{aligned}
	u&(x,t)=\frac{1}{2\pi} \int_\mathbb{R} e^{ikx-\omega(k,t)} \hat{u}_0(k) dk \\
	&- \frac{1}{2\pi} \int_{\partial D_t^+} \frac{e^{ikx-\omega(k,t)}}{\Delta_t(k)} \left[ e^{-ikL} \hat{u}_0(\nu(k))-e^{-i\nu(k)L} \hat{u}_0(k) \right] dk -\frac{1}{2\pi} \int_{\partial D_t^-} \frac{e^{ik(x-L)-\omega(k,t)}}{\Delta_t(k)} \left[\hat{u}_0(\nu(k)) - \hat{u}_0(k) \right] dk\\
	&- \frac{1}{2\pi} \int_{\partial D_t^+} \frac{e^{ikx-\omega(k,t)}}{\Delta_t(k)} \left( 2ik +\frac{C(t)}{B(t)}\right) \left[  e^{-ikL} g_0^b(k,t)-  e^{\frac{C(t)}{B(t)}L} h_0^b(k,t)  \right] dk \\
	&-\frac{1}{2\pi} \int_{\partial D_t^-} \frac{e^{ik(x-L)-\omega(k,t)}}{\Delta_t(k)} \left( 2ik +\frac{C(t)}{B(t)}\right) \left[ g_0^b(k,t)-  e^{i k L} e^{\frac{C(t)}{B(t)}L} h_0^b(k,t)  \right] dk \\
	&+\frac{1}{2\pi} \int_{\partial D_t^+} \frac{e^{ikx}}{\Delta_t(k)} \left[ e^{-ikL} \hat{u}(\nu(k),t)-e^{-i\nu(k)L} \hat{u}(k,t) \right] dk +\frac{1}{2\pi} \int_{\partial D_t^-} \frac{e^{ik(x-L)}}{\Delta_t(k)} \left[ \hat{u}(\nu(k),t)- \hat{u}(k,t) \right] dk,
\end{aligned}
\end{equation}
where $\Delta_t(k)=e^{-ikL}-e^{-i\nu(k)L}, \qquad \nu(k)=-k+i\frac{C(t)}{B(t)}, \quad k\in \mathbb{C}.$

We note that the last two integrals of \eqref{sol-rep-1} vanish. Indeed, let us consider the contribution of the integral
\begin{align*}
\int_{\partial D_t^+} \frac{e^{ikx}}{\Delta_t(k)} e^{-i\nu(k)L} \hat{u}(k,t)  dk = \int_{\partial D_t^+} \frac{e^{ikx+ikL+ \frac{C(t)}{B(t)}L}}{e^{-ikL}-e^{ikL+ \frac{C(t)}{B(t)}L}} \int_0^L e^{-ik\xi} u(\xi,t) d\xi  \,dk 
\end{align*}
The integrand is analytic and bounded in $D_t^+$, since for $k\to\infty$, with $\text{Re}k>0$ its asymptotic behavior is given by
$$ \frac{e^{ikx+ikL+ \frac{C(t)}{B(t)}L}}{e^{-ikL}-e^{ikL+ \frac{C(t)}{B(t)}L}}  \frac{e^{-ikL} u(L,t) - u(0,t)}{ik} \sim \frac{e^{ikx+ \frac{C(t)}{B(t)}L}}{e^{-ikL}-e^{ikL+ \frac{C(t)}{B(t)}L}}  \frac{u(L,t)}{ik}$$
$$  \sim   \frac{1}{e^{-ik(L+x)- \frac{C(t)}{B(t)}L}-e^{ik(L-x)}}  \frac{u(L,t)}{ik} , $$
which is bounded in $D_t^+$ since the denominator displays exponential decay. Indeed, observing  that $k\in \partial D_t^+\Rightarrow \text{Im}k\ge \max \left\{ 0, \frac{C(t)}{B(t)} \right\} $, thus $k\in D_t^+\Rightarrow \text{Im}k > \max \left\{ 0, \frac{C(t)}{B(t)} \right\}$, we obtain $$\text{Re}\left[-ik(L+x)- \frac{C(t)}{B(t)}L \right] = (L+x) \left[ \text{Im}k - \frac{C(t)}{B(t)}\frac{L}{L+x} \right] >  \max \left\{ -L \frac{C(t)}{B(t)} , x \frac{C(t)}{B(t)} \right\} \ge 0.$$ 
We note that $L-x\ge0$, thus $\text{Re}\{i k(L-x)\} =-(L-x) \text{Im}k  \le0$. Similar treatment can be applied to the remaining 3 terms, yielding the desired result. 

Also using the change of variable $k\to\nu(k)$ in the third and fifth integrals one obtains the invariances
\begin{align*}
&\int_{\partial D_t^-} \frac{e^{-\omega(k,t)}}{\Delta_t(k)} dk \ \longrightarrow \ \int_{\partial D_t^+} \frac{e^{-\omega(k,t)}}{\Delta_t(k)} dk \\ 
&\int_{\partial D_t^-} \frac{e^{-\omega(k,t)}}{\Delta_t(k)}\left( 2ik +\frac{C(t)}{B(t)}\right) dk \ \longrightarrow \ -\int_{\partial D_t^+} \frac{e^{-\omega(k,t)}}{\Delta_t(k)} \left( 2ik +\frac{C(t)}{B(t)}\right)dk,
\end{align*}
allowing for the simplified form of the solution that involves integrals only in $\partial D_t^+$:
\begin{equation}\label{sol-rep-2}
\begin{aligned}
	u&(x,t)=\frac{1}{2\pi} \int_\mathbb{R} e^{ikx-\omega(k,t)} \hat{u}_0(k) dk \\
	&- \frac{1}{2\pi} \int_{\partial D_t^+} \frac{e^{-\omega(k,t)}}{\Delta_t(k)} \left[ \left( e^{ik(x-L)} -e^{i\nu(x-L)} \right) \hat{u}_0(\nu)- \left(e^{i k x} -  e^{i \nu x}\right)  e^{-i\nu L} \hat{u}_0(k) \right] dk \\
	&- \frac{1}{2\pi} \int_{\partial D_t^+} \frac{e^{-\omega(k,t)}}{\Delta_t(k)} \left( 2ik +\frac{C(t)}{B(t)}\right) \left[ \left( e^{ik(x-L)}- e^{i\nu(x-L)}\right) g_0^b(k,t)-  \left(e^{i k x} -  e^{i \nu x}\right) e^{\frac{C(t)}{B(t)}L} h_0^b(k,t)  \right] dk.
\end{aligned}
\end{equation}

Introducing the notation $\mu=k-i\frac{C(t)}{2B(t)}$ and recalling  $\nu=-k+i\frac{C(t)}{B(t)}$, rewrites the above solution into the form
\begin{equation}\label{sol-rep-3}
\begin{aligned}
	u&(x,t)=\frac{1}{2\pi} \int_\mathbb{R} e^{ikx-\omega(k,t)} \hat{u}_0(k) dk \\
	&- \frac{ e^{-\frac{C(t)}{2B(t)}x}}{2\pi} \int_{\partial D_t^+} \frac{e^{-\omega(k,t)}}{\sin(\mu L)} \left[\sin[\mu (L-x)]   \hat{u}_0(\nu)+\sin(\mu x)  e^{i\mu L} \hat{u}_0(k) \right] dk \\
	&- \frac{ e^{-\frac{C(t)}{2B(t)}x}}{2\pi} \int_{\partial D_t^+} \frac{e^{-\omega(k,t)}}{\sin(\mu L)} 2i\mu  \left[ \sin[\mu (L-x)] g_0^b(k,t)+  \sin(\mu x) e^{\frac{C(t)}{2B(t)}L} h_0^b(k,t)  \right] dk,
\end{aligned}
\end{equation}
with $\partial D_t^+$ given in Figure \ref{fig:domain}

\begin{remark}\label{rem-1}

In \eqref{sol-rep-2}, can we substitute $g_0^b(k,t)$ with $g_0^b(k,T)$? In other words, what are the conditions that $B(t)$ and $C(t)$ must satisfy so that causality is preserved?

In order to achieve this, we need to show that the contribution of the following integral vanishes

\begin{align*}
	\int_{\partial D_t^+} e^{ikx-\omega(k,t)} \big[g_0^b(k,T) -g_0^b(k,t)\big] dk =  \int_{\partial D_t^+} e^{ikx-\omega(k,t)} \int_t^T e^{\omega(k,s)}b(s) g_0(s) ds dk \\
	=  \int_t^T \int_{\partial D_t^+} e^{ikx+\omega(k,s)-\omega(k,t)} dk \, b(s) g_0(s) ds .
\end{align*}

The inner integral is vanishing if the term $e^{ikx+\omega(k,s)-\omega(k,t)}$ is bounded and analytic when $k\in  D_t^+$. Thus, it is sufficient to show that 
$s>t \Rightarrow \Re \{\omega(k,s)-\omega(k,t)\}<0.$
To investigate this, we will use that  $b(t)\ge b_0>0 \Rightarrow \Big\{ s>t \Rightarrow B(s)> B(t) \Big\}.$ Hence, we would like to prove that

$$\Re \{\omega(k,s)-\omega(k,t)\}<0 \Leftrightarrow \big[B(s)-B(t)\big] 
\left(\Re k^2 + \frac{C(s)-C(t)}{B(s)-B(t)} \Im k\right)<0 \Leftrightarrow 
\Re k^2 \le -\frac{C(s)-C(t)}{B(s)-B(t)} \Im k .$$

Using that $k\in  D_t^+ \Rightarrow \Re \{\omega(k,t)\}\le 0 \Leftrightarrow 
\Re k^2 \le -\frac{C(t)}{B(t)} \Im k$ and $\Im k>0$, a sufficient condition is

$$\frac{C(t)}{B(t)} \ge \frac{C(s)-C(t)}{B(s)-B(t)} \Leftrightarrow \frac{C(s)}{B(s)} \le \frac{C(t)}{B(t)}, \ \text{ if } s>t, $$
namely the function $\frac{C(t)}{B(t)}$ to be non-increasing. Of course, $c(t)=0 \Rightarrow C(t)=0$ or {$C/B=$constant, are special cases} of this condition.

\end{remark}

The following proposition follows from the arguments above in this section.

\begin{proposition}
Let $b=b(t)$ and $c=c(t)$ be two smooth functions on $[0,T]$ (or more globally on $\mathbb{R}$), where $b$ is bounded from below by a strictly positive constant\footnote{In fact, one can rely on the weaker assumption: For each $T>0$, there is $b_0^T>0$ such that $b(t)\ge b_0^T$ for all $t\in [0,T]$.} Let $u=u(x,t)$ satisfy the linear convection-diffusion model \eqref{maineq} on the domain $Q_T= \Omega\times (0,T)$. Assume $u$ is a sufficiently smooth (up to the boundary of $Q_T$) solution of \eqref{maineq}.  Then $u(x,t)$ is given by \eqref{sol-rep-2}. Moreover, $g_0^b(k,t)$ can be substituted with $g_0^b(k,T)$ and $h_0^b(k,t)$ can be substituted with $h_0^b(k,T)$ in \eqref{sol-rep-2} provided $\frac{C(t)}{B(t)}$ is {constant}.
\end{proposition}

\begin{remark}\label{rem-2}[Definition of weak solutions]
Although the solution representation formula is obtained under the assumption that a solution exists and it is smooth, the right hand side of this formula remains to be a valid expression in certain function spaces when initial datum $u_0$ and boundary data $g_0,h_0$ are taken from a rough class of functions (e.g., fractional Sobolev spaces of low order). Therefore, we define also the weak solutions of the associated ibvp through this formula and prove relevant well-posedness results accordingly.  The latter is carried out carefully in Section \ref{wpsection} below.
\end{remark}

\subsection{Reductions}\label{HLsetting}

Letting $L\to\infty$ in \eqref{sol-rep-3}, with the assumption of vanishing solution $u(x,t)$ at $x\to\infty$, namely setting $h_0^b(k,t)=0$  we can retrieve the solution for the case of the half line $\Omega=\mathbb{R}_+$:

\begin{equation}\label{sol-rep-hl}
\begin{aligned}
	u&(x,t)=\frac{1}{2\pi} \int_\mathbb{R} e^{ikx-\omega(k,t)} \hat{u}_0(k) dk \\
	&- \frac{1}{2\pi} \int_{\partial D_t^+} e^{ikx-\omega(k,t)} \left[  \hat{u}_0(\nu)+ \left( 2ik +\frac{C(t)}{B(t)}\right) g_0^b(k,t) \right] dk .
\end{aligned}
\end{equation}

Letting $c(t)=0,\, t>0$, namely $C(t)=0,\,t>0$, in \eqref{sol-rep-3}, we obtain the solution for the heat equation on the finite interval:

\begin{equation}\label{sol-rep-heat-fi}
\begin{aligned}
	u&(x,t)=\frac{1}{2\pi} \int_\mathbb{R} e^{ikx-B(t) k^2} \hat{u}_0(k) dk \\
	&- \frac{1}{2\pi} \int_{\partial D^+} \frac{e^{-B(t) k^2}}{\sin(k L)} \left[\sin[k (L-x)]   \hat{u}_0(-k)+\sin(k x)  e^{ik L} \hat{u}_0(k) \right] dk \\
	&- \frac{1}{2\pi} \int_{\partial D^+} \frac{e^{-B(t) k^2}}{\sin(k L)} 2ik  \left[ \sin[k (L-x)] g_0^b(k,t)+  \sin(k x)  h_0^b(k,t)  \right] dk,
\end{aligned}
\end{equation} where $\partial D^+ = \left\{ k: \Re k^2=0, \ \Im k\ge0 \right\} = \left\{ k: \, k = r+ i |r|, \ r\in\mathbb{R} \right\}$ and 
$$g_0^b(k,t) =\int_0^t b(s) e^{B(s) k^2} g_0(s)ds \ \text{ and } \  h_0^b(k,t) =\int_0^t b(s) e^{B(s) k^2} h_0(s)ds.$$

Letting $L\to\infty$ in \eqref{sol-rep-heat-fi}, with the assumption of vanishing solution $u(x,t)$ at $x\to\infty$, namely setting $h_0^b(k,t)=0$  we can retrieve the solution for the case of the half line:
\begin{equation}\label{sol-rep-heat-hl}
\begin{aligned}
	u&(x,t)=\frac{1}{2\pi} \int_\mathbb{R} e^{ikx-B(t) k^2} \hat{u}_0(k) dk - \frac{1}{2\pi} \int_{\partial D^+} e^{ikx-B(t) k^2} \left[  \hat{u}_0(-k)+  2ik \, g_0^b(k,t) \right] dk .
\end{aligned}
\end{equation}
\subsubsection*{Verification}
Let us consider \eqref{sol-rep-heat-hl}. In view of Remark \ref{rem-1} it is straightforward to check that the PDE is satisfied. Considering the initial condition, and using that $B(0)=0$ and $g_0^b(k,0)=0$, we get:
\begin{align*}
u(x,0)= \frac{1}{2\pi} \int_\mathbb{R} e^{ikx} \hat{u}_0(k) dk - \frac{1}{2\pi} \int_{\partial D^+} e^{ikx}   \hat{u}_0(-k) dk =  \frac{1}{2\pi} \int_\mathbb{R} e^{ikx} \hat{u}_0(k) dk = u_0(x),
\end{align*}
since the integral on $\partial D^+$ is vanishing due to the boundedness and analyticity of the associated integrand.

Considering the boundary condition,
\begin{align*}
u(0,t)= \frac{1}{2\pi} \int_\mathbb{R} e^{-B(t) k^2} \hat{u}_0(k) dk - \frac{1}{2\pi} \int_{\partial D^+} e^{-B(t) k^2}   \hat{u}_0(-k)dk -\frac{1}{2\pi} \int_{\partial D^+} e^{-B(t) k^2}   2ik \, g_0^b(k,t)  dk.
\end{align*}
First, deform the second integral to the real line and make the change of variables $k\to-k$; then it becomes the opposite of the first integral. Second, we make the change of variables $k^2=i\lambda$ in the third integral, thus:
\begin{align*}
u(0,t) &=  \frac{1}{2\pi} \int_\mathbb{R} e^{- i B(t) \lambda}  \,\int_0^t b(s) e^{i B(s) \lambda} g_0(s)ds  d\lambda &=  \frac{1}{2\pi} \int_\mathbb{R} e^{-i B(t) \lambda}  \,\int_0^{B(t)} e^{i r \lambda} g_0(B^{-1}(r))dr  d\lambda \\
&= g_0(B^{-1}(B(t))) = g_0(t),
\end{align*}
where we made the change of variables $r=B(s)$, in the second equality. We note that by definition $B(s)$ is increasing and $dr=b(s)ds$.
\subsection{Illustration of the solution}\label{simul}

In what follows we provide plots of the solutions which display some qualitative difference compared to the IBVPs with constant coefficients. We follow the treatment provided in \cite{AKM24, CFK25} to simplify the associated integrals, while keeping the exponential decay of the integrands.

\begin{figure}[!h]
\centering
\begin{minipage}{0.48\textwidth}
	\centering
	\includegraphics[scale=0.42]{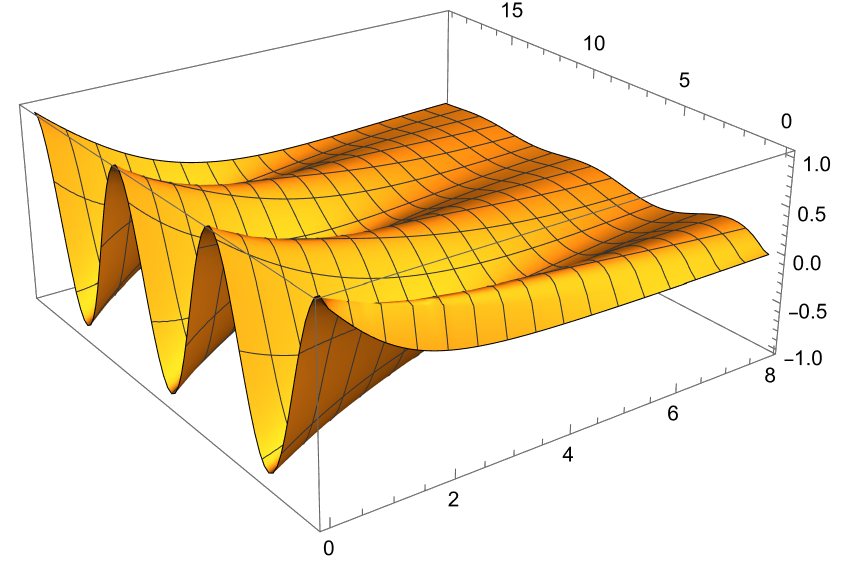}
	\caption{ $u(x,t)$ by \eqref{sol-rep-heat-hl} for $b(t)=\frac{6 \pi }{t+1}+1$  in the domain $(0,8)\times(0,6\pi)$.}
	\label{fig:3D-B}
\end{minipage}\hfill
\begin{minipage}{0.48\textwidth}
	\centering
	\includegraphics[scale=0.4]{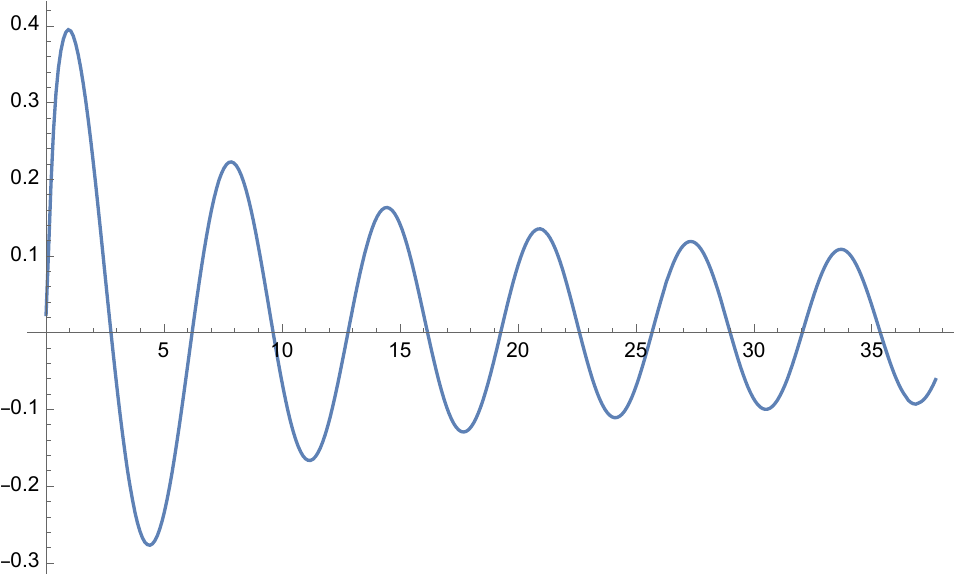}
	\caption{$u(4,t)$ by \eqref{sol-rep-heat-hl} for $b(t)=\frac{6 \pi }{t+1}+1, \ t\in (0,12\pi)$.}
	\label{fig:2D-B}
\end{minipage}
\end{figure}

First, we consider the integral representations \eqref{sol-rep-heat-hl}, with $b(t)=\frac{6 \pi }{t+1}+1$; initial condition $u_0(x)=e^{-x}$; boundary condition $g_0(t)=\cos t$. In Figure \ref{fig:3D-B} we plot $u(x,t)$ in the domain $(0,8)\times(0,6\pi)$. Furthermore, in Figure \ref{fig:2D-B} we plot the $u(4,t), \ t\in(0,12\pi)$, where we observe that the amplitude of the oscillations decreases with $t$; this is consistent with the fact that $b(t)$ is decreasing and $\lim_{t\to\infty} b(t)=1$. 

\begin{figure}[!ht]
\centering
\begin{minipage}{0.48\textwidth}
	\centering
	\includegraphics[scale=0.35]{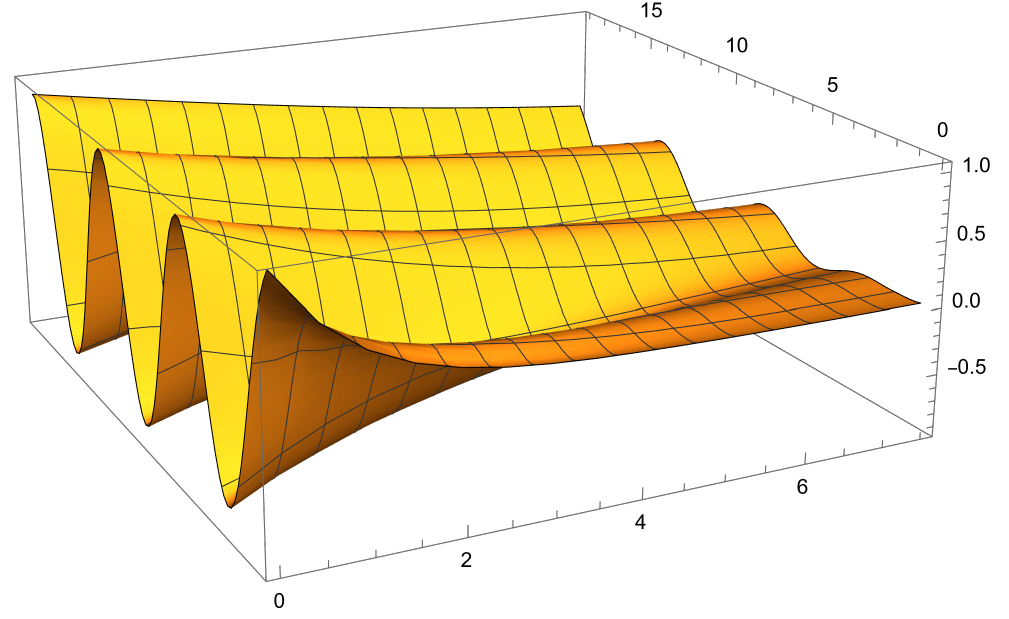}
	\caption{  $u(x,t)$ by \eqref{sol-rep-hl} for $b(t)=1+2t$ and $c(t)=1$ in the domain $(0,8)\times(0,6\pi)$.}
	\label{fig:3D-B-C-1}
\end{minipage}\hfill
\begin{minipage}{0.48\textwidth}
	\centering
	\includegraphics[scale=0.47]{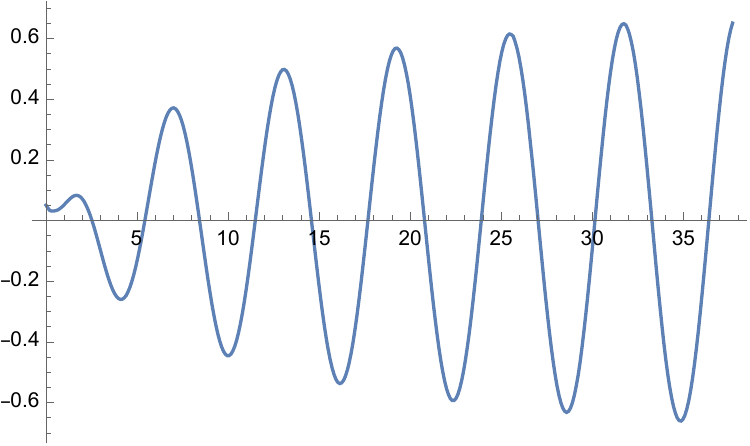}
	\caption{ $u(4,t)$ by \eqref{sol-rep-hl} for $b(t)=1+2t$ and $c(t)=1, \ t\in (0,12\pi)$.}
	\label{fig:2D-B-C-1}
\end{minipage}
\end{figure}

Second, using a different numerical method, we consider the case $b(t)=1+2t$ and $c(t)=1$; initial condition $u_0(x)=e^{-x}$; boundary condition $g_0(t)=\cos t$. In Figure \ref{fig:3D-B-C-1} we plot $u(x,t)$ in the domain $(0,8)\times(0,6\pi)$. Furthermore, in Figure \ref{fig:2D-B-C-1} we plot the $u(4,t), \ t\in(0,12\pi)$, where we observe that the amplitude of the oscillations increases with $t$; this is consistent with the fact that $b(t)$ is increasing, too. 

\begin{figure}[!ht]
\centering
\begin{minipage}{0.48\textwidth}
	\centering
	\includegraphics[scale=0.35]{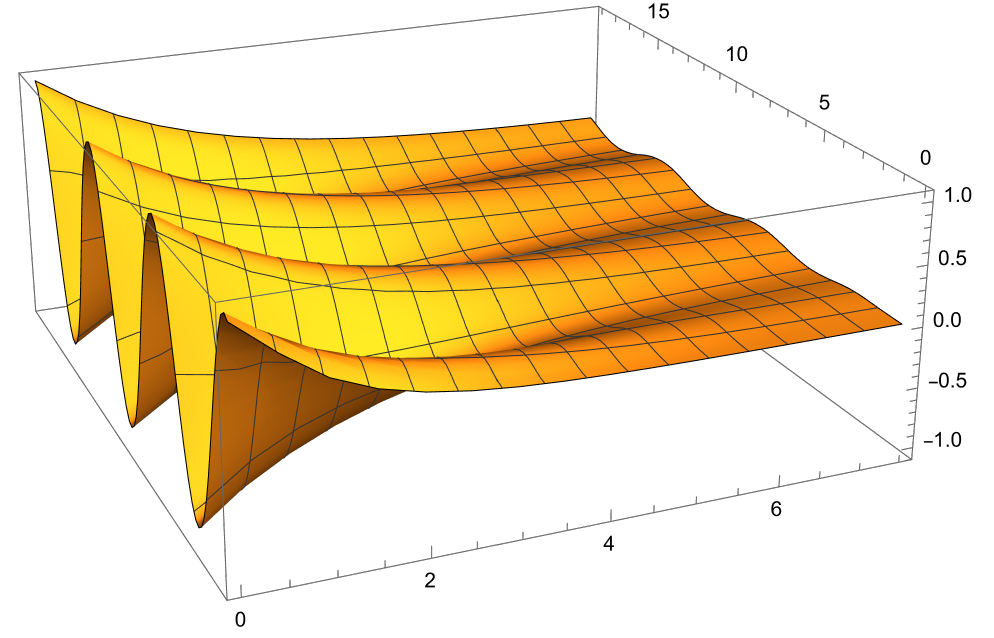}
	\caption{  $u(x,t)$ by \eqref{sol-rep-hl} for $b(t)=1+2t$ and $c(t)=t$ in the domain $(0,8)\times(0,6\pi)$}
	\label{fig:3D-B-C}
\end{minipage}\hfill
\begin{minipage}{0.48\textwidth}
	\centering
	\includegraphics[scale=0.35]{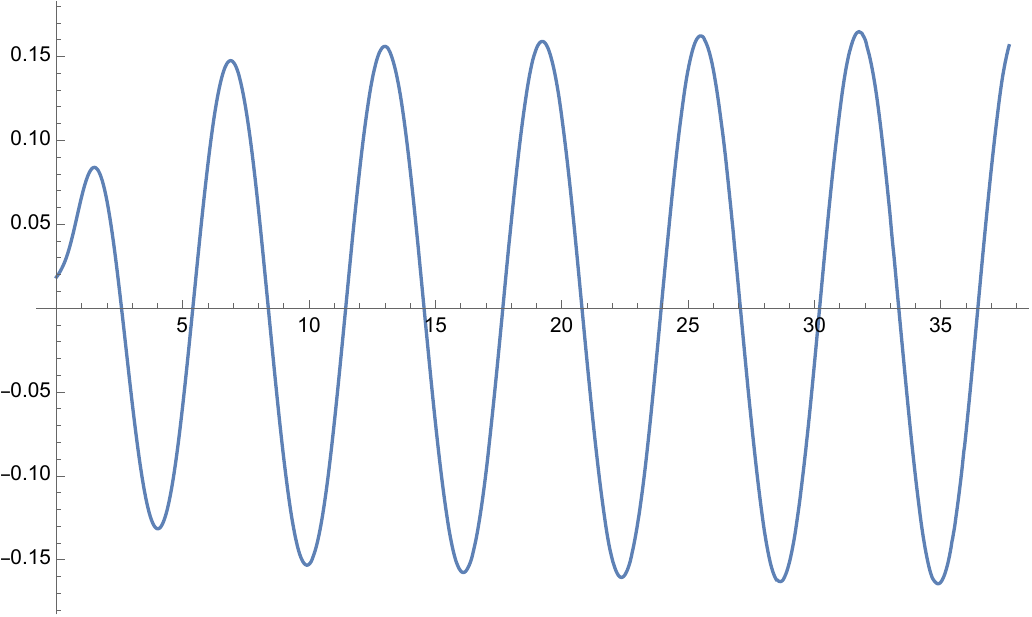}
	\caption{ $u(4,t)$ by \eqref{sol-rep-hl} for $b(t)=1+2t$ and $c(t)=t, \ t\in (0,12\pi)$.}
	\label{fig:2D-B-C}
\end{minipage}
\end{figure}

Third, we consider the case $b(t)=1+2t$ and $c(t)=t$ in the same set-up, with similar observations in Figures \ref{fig:3D-B-C} and \ref{fig:2D-B-C}.

\section{Well-posedness in fractional Sobolev spaces}\label{wpsection}

The well-posedness of the ibvp \eqref{maineq} will be established by using a standard decomposition-reunification algorithm: (Step 1) Extend the initial datum $u_0$ defined on the finite-line to a map $U_0$ defined on the whole-line and solve the associated homogeneous Cauchy problem in \eqref{hom_cauchy} with initial datum $U_0$; (Step 2) Solve a reduced ibvp given by \eqref{reduced_ibvp} with zero initial datum and modified boundary conditions taking into account the boundary traces emerging from the Cauchy problem of Step 1; (Step 3) Finally, recover the solution of the original model through the reunification $u=U|_{(0,L)}+y$, where $U$ solves \eqref{hom_cauchy} with (initial) datum $U_0$ and $y$ solves \eqref{reduced_ibvp} with (boundary) data $(g_0-U|_{x=0}, h_0-U|_{x=L})$.

\subsection{Homogeneous Cauchy problem with extended data}\label{Seccauchy}
We first consider the following homogeneous Cauchy problem on $\mathbb{R}\ni x$:

\begin{equation}\label{hom_cauchy}
U_t=A(t)U,\, t\in (0,T);\quad \quad U(0)=U_0,
\end{equation}
where $U_0\in H^s_x(\mathbb{R})$ with $U_0|_{(0,L)}=u_0$, and $u_0\in H_x^{s}(0,L)$ is the initial datum for the original ibvp in \eqref{maineq}.
To solve \eqref{hom_cauchy}, we take the classical Fourier transform in $x$ and obtain
the infinite family of ODEs:
\begin{equation}
\frac{d}{dt}\hat{U}(k,t)+(k^2b(t)-ikc(t))\hat{U}(k,t)=0, \quad t\in (0,T);
\end{equation} parametrized with $k\in\mathbb{R}$ and for each $k$, the associated ODE is supplemented with the initial condition $\hat{U}(k,0)=\hat{U}_0(k).$  Solving the above initial value problems we obtain
\begin{equation}\label{Uhatkt}
\hat{U}(k,t)=e^{-w(k,t)}\hat{U}_0(k)
\end{equation} for each $k$.  Assuming $b,c\in C^{(n-1)}(0,T)$ for some $n\ge 1$, taking the $m-$th derivative of \eqref{Uhatkt} in $t$, where $1\le m\le n$, in view of Faà di Bruno's formula, we get
\begin{equation}\label{dtmUkt}
\frac{d^m}{dt^m}\hat{U}(k,t)=p_m\left(-\frac{d}{dt}\omega(k,t), -\frac{d^2}{dt^2}\omega(k,t),..., -\frac{d^m}{dt^m}\omega(k,t)\right)e^{-w(k,t)}\hat{U}_0(k),
\end{equation} where $\frac{d^j}{dt^j}\omega(k,t)=k^2b^{(j-1)}(t)-ikc^{(j-1)}(t),\quad 1\le j\le m,$ and $p_m$ is the Bell polynomial, i.e., 
$$p_m(x_1,...,x_m)=m!\sum_{1j_1+2j_2+...+mj_m=m}\prod_{i=1}^m\frac{x_i^{j_i}}{(i!)^{j_i}j_i!}, \quad m\ge 1.$$
For convenience, we set the notation $p_0\equiv 1$. We prove the following spatial regularity results for problem \eqref{hom_cauchy} with the help of the formulation in \eqref{dtmUkt}.
\begin{lemma}[Space estimates and smoothing]\label{hom_lem}
Let $n_1,n_2\ge 1$, $r,s\in\mathbb{R}$, $b\in C^{(n_1-1)}([0,T])$, $c\in C^{(n_2-1)}([0,T])$, $b(t)\ge b_0$ for some $b_0>0$, $u_0\in H^s_x(0,L)$ and $U_0\in H^s_x\mathbb{(R)}$ with that $U_0|_{(0,L)}=u_0$ and\\ $\|U_0\|_{H^s_x(\mathbb{R})}\le 2\|u_0\|_{H^s_x(0,L)}.$ Then, $U\in C([0,T];H^s_x(\mathbb{R}))\cap C^{\min\{n_1,n_2\}}((0,T);H_x^r(\mathbb{R})))$ and it satisfies the estimate 
\begin{equation}\label{space_est_rs1}
	\sup_{t\in [0,T]}\|U(\cdot,t)\|_{H^r_x(\mathbb{R})}\le 2\|u_0\|_{H^s_x(0,L)},\quad \text{ if } r\le s,
\end{equation}
and in addition for each $t>0$ and $0\le m\le \min\{n_1,n_2\}$, there corresponds a constant \\$c_1=c_1({r,s,b_0,m,t,\|b\|_{C^{\min\{n_1,n_2\}-1}([0,T])},\|c\|_{C^{\min\{n_1,n_2\}-1}([0,T])}})>0$ such that
\begin{equation}
	\left\|\frac{\partial^m}{\partial t^m}U(\cdot,t)\right\|_{H^r_x(\mathbb{R})}\le c_1 \|u_0\|_{H^s_x(0,L)}, \quad \text{ if }r>s.
\end{equation}
\end{lemma}
\begin{proof}
We have $\left|\hat{U}(k,t)\right| = e^{-k^2B(t)}|\hat{U}_0(k)|, \quad k\in \mathbb{R},$ and for $0\le m\le \min\{n_1,n_2\}$, 

$$\left|\frac{d^m}{dt^m}\hat{U}(k,t)\right| = e^{-k^2B(t)}\cdot \left|P_{m}^\omega(k,t)\right|\cdot|\hat{U}_0(k)|, \quad k\in \mathbb{R},$$ where

$$P_m^\omega(k,t):= p_m\left(-\frac{d}{dt}\omega(k,t), -\frac{d^2}{dt^2}\omega(k,t),..., -\frac{d^m}{dt^m}\omega(k,t)\right).$$

It follows with any $r>s$ and $t>0$ that

\begin{equation}
	\begin{aligned}
		&\left\|\frac{\partial^m }{\partial t^m}U(\cdot,t)\right\|_{H^r_x(\mathbb{R})}^2 = \int_\mathbb{R}(1+k^2)^r\left|P_m^\omega(k,t)\right|^2\cdot e^{-2k^2B(t)}|\hat{U}_0(k)|^2dk\\
		&=\int_\mathbb{R}(1+k^2)^{r-s}\left|P_m^\omega(k,t)\right|^2e^{-2k^2B(t)}(1+k^2)^s|\hat{U}_0(k)|^2dk\\
		&\le \int_\mathbb{R}(1+k^2)^{r-s}\left|P_m^\omega(k,t)\right|^2e^{-2k^2b_0 t}(1+k^2)^s|\hat{U}_0(k)|^2dk\\
		&\le c_1^2({r,s,b_0,m,t}) \int_\mathbb{R}(1+k^2)^s|\hat{U}_0(k)|^2dk = c_1^2({r,s,b_0,m,t})\|U_0\|_{H^s_x(\mathbb{R})}^2,
	\end{aligned} 
\end{equation} where the last inequality follows from that $(1+k^2)^{r-s}\left|P_m^\omega(k,t)\right|^2e^{-2k^2b_0 t}$ is bounded by a constant $c_1^2$ depending on $r,s\in \mathbb{R}$, $b_0>0$, $m\ge 0$, $t>0$, and bounds of derivatives of $b,c$ to order $\min\{n_1,n_2\}-1$.

On the other hand, if $r\le s$, then we have 

\begin{equation}\begin{aligned}
		\|U(\cdot,t)\|_{H^r_x(\mathbb{R})}^2&\le \|U(\cdot,t)\|_{H^{s}_x(\mathbb{R})}^2=\int_{-\infty}^\infty(1+k^2)^s|\hat{U}(k,t)|^2dk\\
		&\le \int_{-\infty}^\infty(1+k^2)^s|\hat{U}_0(k)|^2dk=\|U_0\|_{H^s_x(\mathbb{R})}^2,\quad \forall t\in [0,T]
	\end{aligned}
\end{equation}
since for each $t\in [0,T]$, we have $|\hat{U}(k,t)|\le |\hat{U}_0(k)|$ for all $k\in\mathbb{R}$.

\end{proof}

\begin{remark}
If $b,c$ are constants, or more generally, $b,c\in C^\infty([0,T])$ in the above lemma, then $U\in C^\infty((0,T);H^r_x(\mathbb{R}))$ for all $r\in \mathbb{R}$.
\end{remark}

\begin{lemma}[Time estimate]
\label{hom_lem_time}
Suppose $b\in C^{n_1}(0,T)\cap C^n([0,T])$, $c\in C^{n_2}(0,T)\cap C^n([0,T])$, $n=\min\{n_1,n_2\}$, and $b(t)\ge b_0$ for some $b_0>0$. Let $u_0\in H^s_x(0,L)$ and $U_0\in H^s_x\mathbb{(R)}$ with that $U_0|_{(0,L)}=u_0$ and $\|U_0\|_{H^s_x(\mathbb{R})}\le 2\|u_0\|_{H^s_x(0,L)}.$ Then, for $-\frac{1}{2}\le s\le {2n-\frac12}$, we have 

$$\sup_{x\in [0,L]}\left\|U(x,\cdot)\right\|_{H^{\frac{2s+1}{4}}_t(0,T)}\le c_2(T,\|b\|_{C^{(n-1)}([0,T])},\|c\|_{C^{(n-1)}([0,T])},n,b_0,s)\cdot \|u_0\|_{H_x^{s}(0,L)},$$
where $c_2$ is a positive constant which remains uniformly bounded as $T\rightarrow 0^+$.

\end{lemma}
\begin{proof}
We take inverse Fourier transform and write $$U(x,t)=\frac{1}{2\pi}\int_{-\infty}^\infty e^{ikx-\omega(k,t)}\hat{U}_0(k)dk,\quad x\in \mathbb{R},\quad t\ge 0.$$ Therefore, 

\begin{equation}
	\begin{aligned}
		\frac{\partial^m}{\partial t^m}U(x,t) = \frac{1}{2\pi}\int_{-\infty}^\infty e^{ikx-\omega(k,t)}P_m^\omega(k,t)\hat{U}_0(k)dk.
	\end{aligned}
\end{equation}
Taking $L_t^2(0,T)$ norm of both sides,
\begin{equation*}
	\begin{aligned}
		&\left\|\frac{\partial^m}{\partial t^m}U(x,\cdot)\right\|_{L^2_t(0,T)}^2 = \int_0^T\left|\frac{1}{2\pi}\int_{-\infty}^\infty e^{ikx-\omega(k,t)}P_m^\omega(k,t)\hat{U_0}(k)dk\right|^2dt\\
		&\le \frac{1}{(2\pi)^2}\int_0^T\left(\int_{-\infty}^\infty e^{-k^2b_0t}\left|P_m^\omega(k,t)\hat{U}_0(k)\right|dk\right)^2dt.
	\end{aligned}
\end{equation*} 

Observe that
\begin{equation}
	\begin{aligned}
		\left|P_m^\omega(k,t)\right|&\le \sum_{\substack{(m_1, \dots, m_m) \in \mathbb{N}_0^m \\
				\sum_{j=1}^m j m_j = m}}
		\frac{m!}{m_1! m_2! \cdots m_m!}
		\prod_{j=1}^m \left(
		\frac{1}{j!} \left( k^2 |b^{(j-1)}(t)| + |k|\cdot |c^{(j-1)}(t)| \right)
		\right)^{m_j}\\
		&\le c_3\left(\|b\|_{C^{(m-1)}([0,T])},\|c\|_{C^{(m-1)}([0,T])},m\right)k^{2m}, \quad |k|>1.\end{aligned}
\end{equation} 

It follows from splitting integration domain in $k$ into $|k|\le 1$ and $|k|>1$ that
\begin{equation*}
	\begin{aligned}
		\left\|\frac{\partial^m}{\partial t^m}U(x,\cdot)\right\|_{L^2_t(0,T)}^2 &\le \frac{c_3^2}{(2\pi)^2b_0}
		\int_0^{b_0T}\left(\int_1^\infty e^{-k^2t'}k^{2m}\left(|\hat{U}_0(k)|+|\hat{U}_0(-k)|\right)dk\right)^2dt'\\
		&+c_4^2\left(\|b\|_{C^{(m-1)}([0,T])},\|c\|_{C^{(m-1)}([0,T])},m,b_0,s\right)T\|U_0\|_{H_x^{s}(\mathbb{R})}^2, \quad s\in\mathbb{R}.
	\end{aligned}
\end{equation*}

Hence, using arguments relying on the boundedness of the Laplace transform (see e.g., \cite[Theorem 3.1]{HMY19} for the case $m=0$), the following estimate follows $$\left\|\frac{\partial^m}{\partial t^m}U(x,\cdot)\right\|_{L^2_t(0,T)}\lesssim (c_3 +\sqrt{T}c_4)\|U_0\|_{H_x^{s}(\mathbb{R})}, \forall s\ge 2m-\frac{1}{2}.$$

Summing up the estimates from $m=0$ to $m=n$, and interpolating we get

$$\left\|U(x,\cdot)\right\|_{H^{\frac{2s+1}{4}}_t(0,T)}\lesssim c_2\|U_0\|_{H_x^{s}(\mathbb{R})}, \quad -\frac{1}{2}\le s \le  {2n-\frac12},$$

where $c_2= c_2(T,\|b\|_{C^{(n-1)}([0,T])},\|c\|_{C^{(n-1)}([0,T])},n,b_0,s)$ is a positive constant which remains uniformly bounded as $T\rightarrow 0^+$.
\end{proof}

\begin{remark}
The above lemma extends \cite[Theorem 3.1]{HMY19} to time dependent coefficients case, as well as to the range $s>\frac32$ when $b$ and $c$ are sufficiently smooth. Moreover, one can improve the boundary regularity estimate of \cite[Theorem 3.1]{HMY19} further through some gain due to smoothing property of the heat semigroup. For instance, we have

\begin{equation}
	\begin{aligned}\label{boostest1}
		&\left\|U(x,\cdot)\right\|_{L^2_t(0,T)} = \left(\int_0^T\left|\frac{1}{2\pi}\int_{-\infty}^\infty e^{ikx-\omega(k,t)}\hat{U}_0(k)dk\right|^2dt\right)^{\frac{1}{2}}\\
		&\le \frac{1}{(2\pi)}\int_0^T\left(\int_{-\infty}^\infty e^{-2k^2b_0t}\left|\hat{U}_0(k)\right|^2dk\right)^\frac{1}{2}dt\\
		& \lesssim  \int_0^T\left(\int_{|k|\le 1} e^{-2k^2b_0t}\left|\hat{U}_0(k)\right|^2dk\right)^\frac{1}{2}dt + \int_0^T\left(\int_{|k|>1} e^{-2k^2b_0t}\left|\hat{U}_0(k)\right|^2dk\right)^\frac{1}{2}dt.
	\end{aligned}
\end{equation}

The first integral above in the last line of \eqref{boostest1}, corresponding to the low frequency modes, can be estimated as 

\begin{equation}
	\begin{aligned}
		\int_0^T&\left(\int_{|k|\le 1} e^{-2k^2b_0t}\left|\hat{U}_0(k)\right|^2dk\right)^\frac{1}{2}dt \\
		& = \int_0^T\left(\int_{|k|\le 1} (1+k^2)^{-s}e^{-2k^2b_0t}(1+k^2)^s\left|\hat{U}_0(k)\right|^2dk\right)^\frac{1}{2}dt \lesssim_{s,T} \|U_0\|_{H_x^s(\mathbb{R})} 
	\end{aligned} 
\end{equation} for any $s\in\mathbb{R}$, while the second integral in the last line of \eqref{boostest1} can be estimated as

\begin{equation}
	\begin{aligned}
		\int_0^T &\left(\int_{|k|>1} e^{-2k^2b_0t}\left|\hat{U}_0(k)\right|^2dk\right)^\frac{1}{2}dt \le \int_0^T\left(\int_{|k|>1} 2k^{2\rho}e^{-2k^2b_0t} (1+k^2)^{-\rho}\left|\hat{U}_0(k)\right|^2dk\right)^\frac{1}{2}dt\\
		&\lesssim_{b_0,\rho}\left(\int_0^Tt^{-\rho/2}dt\right)\cdot \|U_0\|_{H_x^{-\rho}(\mathbb{R})}\lesssim_{b_0,\rho,T}\|U_0\|_{H_x^{-\rho}(\mathbb{R})},
	\end{aligned}
\end{equation} provided $\rho<2$. Therefore, one has $\left\|U(x,\cdot)\right\|_{L^2_t(0,T)}  \lesssim_{b_0,\epsilon,T} \|U_0\|_{H_x^{-2+\epsilon}(\mathbb{R})}$ for any (arbitrarily small) $\epsilon>0$. Comparing this with the result of \cite[Theorem 3.1]{HMY19}, where the right hand side norm is $H^{-\frac12}_x(\mathbb{R})$, one gains $-\frac12+2-\epsilon=\frac{3}{2}-\epsilon$ derivatives for $L^2_t(0,T)$ estimate of the Dirichlet traces.
\end{remark}

\subsection{Reduced initial boundary value problem}\label{Secredibvp}
In this section, we study the following model (with zero initial datum and modified Dirichlet traces):
\begin{equation}\label{reduced_ibvp}
\begin{aligned}
	y_t &= A(t)y; \quad t\in (0,T');\\
	y(0) &= 0;\quad
	\gamma_0 y=(g_0-U(0,\cdot),h_0-U(L,\cdot)),
\end{aligned}
\end{equation} where $T'>T$, $U$ is the solution of the associated homogeneous Cauchy problem that was studied in the previous section on the time interval $[0,T]$ with data $u_0\in H_x^s(0,L)$. To simplify the notation, let us set $g:=g_0-U(0,\cdot)$ and $h:=h_0-U(L,\cdot)$, where $g_0,h_0\in H_t^{\frac{2s+1}{4}}(0,T)$.  Recall that $U(0,\cdot), U(L,\cdot)\in H_t^{\frac{2s+1}{4}}(0,T)$ when $u_0\in H^s_x(0,L)$ (see Lemma \ref{hom_lem_time}). Here, we apply a particular bounded extension operator to $g_0-U(0,\cdot)$ and $h_0-U(L,\cdot)$, extending them from $[0,T]$ to $\mathbb{R}_t$ so that $\text{sup} g, \text{sup} h\subseteq [0,T')$ and also $\label{bd_ext_op}\|g\|_{H_t^{\frac{2s+1}{4}}(0,T')}\le \|g\|_{H_t^{\frac{2s+1}{4}}(\mathbb{R})}\lesssim \|g\|_{H_t^{\frac{2s+1}{4}}(0,T)}$ and similarly 
$\label{bd_ext_op2}\|h\|_{H_t^{\frac{2s+1}{4}}(0,T')}\le \|h\|_{H_t^{\frac{2s+1}{4}}(\mathbb{R})}\lesssim \|h\|_{H_t^{\frac{2s+1}{4}}(0,T)},$ where $s\ge -\frac12$.  In order to guarantee existence of such extensions without distorting regularity, it is necessary to assume that initial-boundary data pairs $(u_0,g_0)$ and $(u_0,h_0)$ satisfy the necessary compatibility conditions (these are local compatibility conditions at space-time corner points $(x,t)=(0,0)$ and $(x,t)=(L,0)$ when $s>\frac12$, $s\neq \frac{2n+1}{2}$, $n\in\mathbb{N}$).  We also assume in this section that $b,c\in C^{\lceil s\rceil}([0,T])$, $b\ge b_0>0$, in particular we assume these properties on $[0,T']$.

We will assume without loss of generality $h\equiv 0$ because the general case can be considered by splitting \eqref{reduced_ibvp} into two problems below:
\begin{equation}\label{reduced_ibvp_q}
q_t = A(t)q; \quad t\in (0,T');\quad \quad  q(0) = 0;\quad \gamma_0 q=(g,0),
\end{equation}
\begin{equation}\label{reduced_ibvp_z}
z_t = A(t)z; \quad t\in (0,T');\quad \quad z(0) = 0;\quad
\gamma_0 z=(0,h).
\end{equation} Furthermore the second problem above can be transformed into one of type \eqref{reduced_ibvp_q} (with $A(t)$ replaced with $\tilde{A}(t) =b(t)\partial_x^2-c(t)\partial_x$  and $g$ replaced with $h$) via the change of variable $$v(x,t):=z(L-x,t),\quad x\in [0,L],\, t\in [0,T').$$

We have the integral representation for the solution of \eqref{reduced_ibvp_q} (with $h\equiv 0$) given by the formula

\begin{equation*}
q(x,t)=- \frac{1}{2\pi} \int_{\partial D_t^+} \frac{e^{-\omega(k,t)}}{\Delta_t(k)} \left( 2ik +\frac{C(t)}{B(t)}\right)  \left( e^{ik(x-L)}- e^{i\nu(k)(x-L)}\right) g^b(k,T')  dk.
\end{equation*}

assuming $\frac{C(t)}{B(t)}$ is {constant}. We parametrize $\partial D_t^+$ with a real parameter via
$k(\lambda) = \lambda+i\left(\frac{C(t)}{2B(t)}+\sqrt{\lambda^2+\frac{C^2(t)}{4B^2(t)}}\right), \quad \lambda\in (-\infty,\infty).$
Therefore, $k$ is differentiable with respect to $\lambda$ on $(-\infty,0)\cup (0,\infty)$, indeed on $\mathbb{R}$ if $C(t)\neq 0$. Note when $C(t)\equiv 0$, $\partial D_t^+$ is union of two rays $k=\lambda+i\lambda$, $\lambda\in [0,\infty)$ and $k=\lambda-i\lambda$, $\lambda\in (-\infty,0]$, with a corner at zero, where $k$ is not differentiable with respect to $\lambda$. Therefore, it is convenient to split the integration region into $I_1=(0,\infty)$ and $I_2=(-\infty,0)$. Thus, we can write the solution formula as follows
\begin{equation*}\label{sol-rep-reduced}
\begin{aligned}
	q(x,t) &=
	- \frac{1}{2\pi} \sum_{j=1}^2 \int_{I_j}  \frac{e^{-\omega(k(\lambda),t)}}{\Delta_t(k(\lambda))} \left( 2ik(\lambda) +\frac{C(t)}{B(t)}\right)  \left( e^{ik(\lambda)(x-L)}- e^{i\nu(k(\lambda))(x-L)}\right) g^b(k(\lambda),T') \frac{dk(\lambda)}{d\lambda}   d\lambda\\
	&=:q_1(x,t)+q_2(x,t),
\end{aligned}
\end{equation*} where $q_1,q_2$ are the integrals associated with $I_1$ and $I_2$, respectively.  We make some observations. At first, if $C(t)\neq 0$, we have 
$
\frac{dk(\lambda)}{d\lambda} = 1+i\frac{\lambda}{\sqrt{\lambda^2+\frac{C^2(t)}{4B^2(t)}}}, \quad \lambda\in \mathbb{R},
$ and if $C(t)=0$, we have $\frac{dk(\lambda)}{d\lambda}=(1+\text{sign}(\lambda)i), \lambda\neq 0.$  It follows that for all $\lambda\in I_1\cup I_2$, we have $\left|\frac{dk(\lambda)}{d\lambda}\right|\le 2.$

Refraining from extension in $t$ beyond $T'$, one gets the estimate
\begin{equation}
\begin{aligned}
	\left|g^b(k(\lambda),T')\right|\le\left|\int_0^{T'}e^{\omega(k(\lambda),s)}b(s)g(s)ds\right|\le \int_0^{T'}|b(s)g(s)|ds \le \|b\|_{L^2_t(0,T')}\cdot\|g\|_{L^2_t(0,T')}.
\end{aligned}
\end{equation}

Since $\frac{C(t)}{B(t)}$ is non-increasing, in view of l'Hôpital and $b(t)\ge b_0>0$, a crude estimate yields
$$\left|2ik(\lambda) +\frac{C(t)}{B(t)}\right|\le 2|k(\lambda)|+\frac{|c(0)|}{b_0}\le \max\left\{4,\frac{3|c(0)|}{b_0}\right\}\left(|\lambda|+1\right),\quad  \forall \lambda\in \mathbb{R}.$$

Recall that we have $\text{Re}(\omega(k,t))=0$ on $\partial D_t^+$, it follows that $|e^{-\omega(k(\lambda),t)}|=1$.  Now, we want to estimate $L^2_x(0,L)$ norm of $q_1(\cdot,t)$. To this end, we split the integration region of $q_1$ into two parts, $0<\lambda\le 1$ and $\lambda>1$. Let us write $q_1=q_{1,1}+q_{1,2}$, where $q_{1,1}$ and $q_{1,2}$ are the integrals over the aforementioned regions, respectively. Then, we have, upon using previous estimates and in addition changing the spatial variable via $x'=L-x$,
\begin{equation}\label{q11smallambda}
\begin{aligned}
	\|&q_{1,1}(\cdot,t)\|_{L^2_x(0,L)}^2\\
	&= \frac{1}{2\pi}\int_0^L\left|\int_{0}^1  \frac{e^{-\omega(k(\lambda),t)}}{\Delta_t(k(\lambda))} \left( 2ik(\lambda) +\frac{C(t)}{B(t)}\right)  \left( e^{ik(\lambda)(x-L)}- e^{i\nu(k(\lambda))(x-L)}\right) g^b(k(\lambda),T') \frac{dk(\lambda)}{d\lambda}   d\lambda\right|^2dx\\
	&\lesssim \int_0^L\int_0^1|g^b(k(\lambda),T')|^2d\lambda dx'\lesssim \|g\|_{L^2(0,T')}^2\lesssim \|g\|_{L^2(0,T)}^2,
\end{aligned}
\end{equation} where we also used
$\displaystyle \frac{ \left|e^{-ik(\lambda)x'}- e^{-i\nu(k(\lambda))x'}\right|}{\left|e^{-ik(\lambda)L}-e^{-i\nu(k(\lambda))L}\right|}\le 1
$ (uniformly bounded) and Cauchy-Schwarz inequality.

Next, we want to estimate $q_{1,2}$, i.e., $\lambda>1$ regime.  Using definition of $g^b$ and the fact that it is compactly supported, we have \begin{equation}
\begin{aligned}
	g^b(k(\lambda),T')=\int_0^{T'}e^{w(k(\lambda),s)}b(s)g(s)ds=\int_{-\infty}^\infty e^{w(k(\lambda),s)}b(s)g(s)ds\\
	=\int_{-\infty}^\infty e^{2iB(s)\lambda\sqrt{\lambda^2+\frac{C^2(s)}{4B^2(s)}}}b(s)g(s)ds.
\end{aligned}
\end{equation}

\textbf{Case 1:} First, we consider the case $c(t)\equiv 0$. In this case,
\begin{equation}
\begin{aligned}\label{Ghatfirst}
	g^b(k(\lambda),T')&=\int_{0}^{T'} e^{2iB(s)\lambda^2}b(s)g(s)ds=\int_{0}^{B(T')}e^{2i\tau\lambda^2}G(\tau)d\tau=\hat{G}(-2\lambda^2),
\end{aligned}
\end{equation}where $G(\tau):=g(\tilde{B}^{-1}(\tau))$, where $\tilde{B}(t)$ is a smooth monotone increasing extension of $B(t)$ to $\mathbb{R}$.

Recall that the $n$th derivative of a linear composition $X_\phi(g):=g\circ \phi$ is given by the Faà di Bruno's formula:
\[
\frac{d^n}{ds^n} g(\phi(s)) = 
\sum \frac{n!}{m_1! 1!^{m_1} m_2! 2!^{m_2} \cdots m_n! n!^{m_n}} \cdot g^{(k)}(\phi(x)) \cdot 
\prod_{j=1}^{n} \left( \phi^{(j)}(s) \right)^{m_j}
\]
where $\phi=B^{-1}$ and the sum is over all sequences of non-negative integers $(m_1, m_2, \dots, m_n)$ satisfying:
$
\sum_{j=1}^{n} m_j = k, \quad \sum_{j=1}^{n} j m_j = n.
$

It follows from this formula that 
\begin{align*}
\|X_\phi(g)\|_{L^2(\mathbb{R})}^2=\int_0^{B(T')}|X_\phi(g)|^2ds =\int_0^{B(T')}|g(B^{-1}(s))|^2ds=\int_0^{T'}|g(t)|^2b(s)ds\le \|b\|_{L^\infty(0,T')}\|g\|_{L^2(0,T')}^2
\end{align*}
and for any $n\ge 1$,
\begin{equation}
\begin{aligned}
	&\left\|\frac{d^n}{ds^n}X_\phi(g)\right\|_{L^2(\mathbb{R})}^2=\int_0^{B(T')}\left|\frac{d^n}{ds^n}X_\phi(g)\right|^2ds =\int_0^{B(T')}\left|\frac{d^n}{ds^n}g(B^{-1}(s))\right|^2ds\\
	&\le c_0(n, \|B^{-1}\|_{C^n([0,T'])}, \|b\|_\infty)\sum_{j=1}^n\|g^{(j)}\|_{L^2(0,T')}^2 \equiv c_0(n, \|B^{-1}\|_{C^n([0,T'])}, \|b\|_\infty)\|g\|_{H^n(0,T')}^2.
\end{aligned}
\end{equation}
Since, $X_\phi$ is linear, we can interpolate, and therefore more general, i.e., fractional cases can be included in the above estimate:
\begin{equation}
\|G\|_{H^\alpha(\mathbb{R})} = \left\|X_\phi(g)\right\|_{H^\alpha(\mathbb{R)}}\lesssim c_0(\lceil \alpha\rceil, \|B^{-1}\|_{C^{\lceil \alpha\rceil}([0,T'])}, \|b\|_\infty)\|g\|_{H^\alpha(0,T')}, \quad \alpha\ge 0.
\end{equation}

Observe that since $\lambda>1$ and $C(t)\equiv 0$, we have 
$$\left|\frac{1}{\Delta_t(k(\lambda))}\right|=\left|\frac{1}{e^{-ik(\lambda)L}-e^{-i\nu(k(\lambda))L}}\right|\le \frac{1}{e^{\lambda L}-e^{-\lambda L}}\le \frac{e^{-\lambda L}}{1-e^{-2\lambda L}}$$
and $\left|e^{ik(\lambda)(x-L)}-e^{-i\nu(k(\lambda))(x-L)}\right|\le e^{\lambda x'}-e^{-\lambda x'}$ with $x'=L-x$.  Using above inequalities, we get 
\begin{equation}\label{q12largelambda}
\begin{aligned}
	\|&q_{1,2}(\cdot,t)\|_{L^2_x(0,L)}^2\\
	&= \frac{1}{2\pi}\int_0^L\left|\int_{1}^\infty  \frac{e^{-\omega(k(\lambda),t)}}{\Delta_t(k(\lambda))} \, 2ik(\lambda)\,  \left( e^{-ik(\lambda)x'}- e^{i\nu(k(\lambda))x'}\right) g^b(k(\lambda),T') \frac{dk(\lambda)}{d\lambda}   d\lambda\right|^2dx\\
	&\lesssim \int_0^L\left(\int_1^\infty {e^{-\lambda x'}}\lambda|\hat{G}(-2\lambda^2)|d\lambda\right)^2 dx'\lesssim \int_1^\infty \lambda^2|\hat{G}(-2\lambda^2)|^2d\lambda \lesssim \|G\|_{H_t^{\frac{1}{4}}(\mathbb{R})}^2 
	\lesssim \|g\|_{H_t^{\frac{1}{4}}({0,T})}^2.
\end{aligned}
\end{equation}

Above, we established $L_x^2$-level estimate for the solution. Next, we will prove an $H_x^n$ level estimate. To this end, we differentiate the solution formula $n$-times with respect to $x$. We have for instance, with $c(t)\equiv 0$,
\begin{equation*}\label{q1nthder}
\begin{aligned}
	\partial_x^nq_1(x,t) &=
	\frac{2}{\pi} \int_{I_1}  \frac{e^{-\omega(k(\lambda),t)}}{\Delta_t(k(\lambda))} \lambda  \left( (ik(\lambda))^ne^{ik(\lambda)(x-L)}- (i\nu(k(\lambda)))^ne^{i\nu(k(\lambda))(x-L)}\right) g^b(k(\lambda),T)  d\lambda,
\end{aligned}
\end{equation*}
where $k(\lambda)=(1+i)\lambda=-\nu(k(\lambda))$. It is not difficult to show that, uniformly in $x\in[0,L]$,
\begin{equation}
\left| \frac{(ik(\lambda))^ne^{ik(\lambda)(x-L)}- (i\nu(k(\lambda)))^ne^{i\nu(k(\lambda))(x-L)}}{\Delta_t(k(\lambda))}\right|\lesssim \begin{cases}\lambda^{n-1} \text{ as }\lambda\rightarrow 0^+ \text{ for }n \text{ odd};\\
	\lambda^{n} \text{ as }\lambda\rightarrow 0^+ \text{ for }n \text{ even};\\
	\lambda^{n} \text{ as }\lambda\rightarrow \infty \text{ for }n \text{ even or odd}.
\end{cases}
\end{equation}

Therefore, \eqref{q11smallambda} is similar for the $L_x^2$-norm of the $n$-th  spatial derivative.
Also, when estimating the $L_x^2$-norm of $n$-th  spatial derivative, the integral $\int_1^\infty \lambda^2|\hat{G}(-2\lambda^2)|^2d\lambda$ in \eqref{q12largelambda} should be replaced with \\$\int_1^\infty \lambda^{2+2n}|\hat{G}(-2\lambda^2)|^2d\lambda$ which can be bounded by $\|g\|_{H_t^{\frac{2n+1}{4}}({0,T})}^2$. Finally, by summing estimates for the first $n$ derivatives and interpolating, applying also the same arguments to $q_2$, we establish the estimate
$               \|q(\cdot,t)\|_{H^s_x(0,L)} \lesssim \|g\|_{H_t^{\frac{2s+1}{4}}({0,T})}, s\ge 0,
$ where the constant depends only on the fixed parameters of the problem, and $b$ is assumed to be from the class of functions $C^{(\lceil s\rceil)}([0,T])$.
So, we just proved the following lemma in view of above estimates and also Lemma \ref{hom_lem_time}.
\begin{lemma}[Ibvp space estimate with zero drift]
\label{zero_lem}
Let $T'>T>0$, $s\ge 0$, $s\neq \frac{2n+1}{2}$ for $n\in\mathbb{N}$, $c(t)\equiv 0$, $b\in C^{(\lceil s\rceil)}([0,T'])$, $b(t)\ge b_0$ for some $b_0>0$, $g_0,h_0\in H_t^{\frac{2s+1}{4}}(0,T)$, $u_0\in H_x^{s}(0,L)$ so that $U(0,\cdot), U(L,\cdot)\in H_t^{\frac{2s+1}{4}}(0,T)$.  Assume that $(u_0,g_0)$ and $(u_0,h_0)$ satisfy the necessary compatibility conditions at the space-time corner points $(x,t)=(0,0)$ and $(x,t)=(L,0)$, respectively. Then, the solution of \eqref{reduced_ibvp} restricted to $(0,T)$ (still denoted same) satisfies the space estimate
\begin{equation}\label{space_est_y}
	\sup_{t\in [0,T]}\|y(\cdot,t)\|_{H^s_x(0,L)}\lesssim \|g_0\|_{H_t^{\frac{2s+1}{4}}({0,T})} + \|h_0\|_{H_t^{\frac{2s+1}{4}}({0,T})} + c_T\|u_0\|_{H_x^s(0,L)}, \quad s\ge 0,
\end{equation} where the constant of the inequality depends on $s,b_0,$ and the bounds of first $\lceil s\rceil$ derivatives of $b$, and the bound of the fixed extension operator leading to estimates \eqref{bd_ext_op}-\eqref{bd_ext_op2}, moreover $c_T$ remains uniformly bounded as $T\rightarrow 0^+$.
\end{lemma}
\begin{remark}
We note that in the case of non-zero drift, since the quantity $C^2/4B^2$ is nonzero, the last identity at the right hand side of \eqref{Ghatfirst} is no longer valid.  Therefore, is is not clear how to relate the associated spectral function to a fractional Sobolev norm as in \eqref{q12largelambda}.  Thus, we will follow an indirect approach for the nonzero drift case (Case 2 below).
\end{remark}
\textbf{Case 2:} Now, we consider the case $c(t)\not\equiv  0$.  To this end, we first consider the following non-homogeneous Cauchy problem ($c(t)\equiv 0$) with zero initial datum:
\begin{equation}\label{non_hom_cauchy}
W_t=b(t)W_{xx} + F_x(x,t),\, t\in (0,T);\quad \quad W(0)=0,
\end{equation} where $F$ is a fixed function.
Using Duhamel's principle, the solution of above Cauchy problem can be written as 
$
W(x,t) = \int_0^tS(t-t')[F_x(\cdot,t')]dt',
$ where $S(t)[\varphi]$ denotes the solution of homogeneous Cauchy problem \eqref{hom_cauchy} with $c(t)\equiv 0$ and initial datum $\varphi$.

We will first prove an estimate for the semigroup $S(t)$.  We have
$\mathcal{F}(S(t)\varphi)(k) = e^{-k^2B(t)}\hat{\varphi}(k).$ Therefore, $\mathcal{F}(\partial_x(S(t)\varphi)) = ike^{-k^2B(t)}\hat{\varphi}(k).$ It follows that
\begin{equation}
\begin{aligned}
	\|\partial_x(S(t)\varphi)\|_{H_x^s(\mathbb{R})}^2 &= \int_{-\infty}^\infty (1+k^2)^{s}k^2e^{-2k^2B(t)}|\hat{\varphi}(k)|^2dk\\
	&\le \frac{e^{-1}}{{2B(t)}}\int_{-\infty}^\infty (1+k^2)^{s}|\hat{\varphi}(k)|^2dk \le \frac{\|\varphi\|_{H_x^s(\mathbb{R})}^2}{{2eb_0t}}.
\end{aligned}
\end{equation} since the maximum value of the function $k^2e^{-2k^2B(t)}$ with respect to $k$ is $\frac{e^{-1}}{2B(t)}$.

Using this estimate, we obtain
\begin{equation}
\begin{aligned}
	\|W(\cdot,t)\|_{H^s_x(\mathbb{R})}&\le \int_0^t\left\|S(t-t')[F_x(\cdot,t')]\right\|_{H_x^s(\mathbb{R})}dt'=\int_0^t\left\|\partial_x(S(t-t')[F(\cdot,t')])\right\|_{H_x^s(\mathbb{R})}dt'\\
	&\le \frac{e^{-1/2}}{\sqrt{2b_0}}\int_0^t\frac{1}{\sqrt{t-t'}} \|F(\cdot,t')\|_{H^s_x(\mathbb{R})}dt'\le\frac{2e^{-1/2}\sqrt{T}}{\sqrt{2b_0}} \|F\|_{C([0,T];H^s_x(\mathbb{R}))},
\end{aligned}
\end{equation}which implies $\|W\|_{C([0,T];H^s_x(\mathbb{R}))}\lesssim \sqrt{T}\|F\|_{C([0,T];H^s_x(\mathbb{R}))}, \quad s\in \mathbb{R}.$
\begin{remark}
Notice that we have $F_x$ at the right hand side of the equation instead of $F$. Therefore, the above estimate improves \cite[Theorem 3.2]{HMY19} by one unit in $s$, also extends it to time dependent coefficient case.  This becomes possible due to utilization of the smoothing property of the heat kernel.
\end{remark}

Next, we will prove time trace estimates for the non-homogeneous Cauchy problem but again utilizing the smoothing property of the heat semigroup associated with the operator $b(t)\partial_x^2$.  We start with the $L_t^2(0,T)$ estimate:

\begin{equation}\label{WL2t}
\begin{aligned}
	\|W(x,\cdot)\|_{L^2_t(0,T)} &= \left(\int_0^T|W(x,t)|^2dt\right)^\frac12 = \left(\int_0^T\left|\gamma_x\int_0^tS(t-t')[F_x(\cdot,t')]dt'\right|^2dt\right)^{\frac12}\\
	& \le \int_0^T\left(\int_0^t\left|\gamma_xS(t-t')[F_x(\cdot,t')]\right|^2dt'\right)^{\frac12}dt \\
	&=\int_0^T\left(\int_0^t\left|\gamma_x\partial_x(S(t-t')[F(\cdot,t')])\right|^2dt'\right)^{\frac12}dt.
\end{aligned}
\end{equation}

We consider the trace at $x$ of the spatial derivative of solution of homogeneous Cauchy problem:
$\gamma_x\partial_x(S(t)\varphi) = \frac{1}{2\pi}\int_{-\infty}^\infty ike^{ikx-k^2B(t)}\hat{\varphi}(k) dk.$ We have
\begin{equation}
\begin{aligned}
	\left(\int_0^t\left|\gamma_x\partial_x(S(t)\varphi)\right|^2dt\right)^{\frac12} &= \left(\int_0^t \left|\frac{1}{2\pi}\int_{-\infty}^\infty ike^{ikx-k^2B(t)}\hat{\varphi}(k) dk\right|^2dt\right)^{\frac12}\\
	& \le \frac{1}{2\pi} \int_0^t \left(\int_{-\infty}^\infty k^2e^{-2k^2B(t)}|\hat{\varphi}(k)|^2dk\right)^{\frac12}dt\\
	& \le \frac{1}{2\pi} \int_0^t \left(\int_{-\infty}^\infty k^2(1+k^2)^{\rho} e^{-2k^2B(t)}(1+k^2)^{-\rho}|\hat{\varphi}(k)|^2dk\right)^{\frac12}dt\\
	&\le \frac{1}{2\pi} \int_0^t \left[\frac{(1+\rho)^{\frac{1+\rho}{2}}}{(2eB(t))^{\frac{1+\rho}{2}}}+\frac{2^\rho}{\sqrt{2eB(t)}}\right]\left(\int_{-\infty}^\infty(1+k^2)^{-\rho} |\hat{\varphi}(k)|^2dk\right)^{\frac12}dt\\
	&\lesssim ({t^{\frac{(1-\rho)}{2}}+t^\frac12)\|\varphi\|_{H^{-{\rho}}_x(\mathbb{R})}},
\end{aligned}
\end{equation} where $\rho=3-\epsilon$ with $\epsilon>0$ arbitrarily small. Using the last estimate in \eqref{WL2t}, we get

$$\|W(x,\cdot)\|_{L^2_t(0,T)}\lesssim  (T^{\frac{3-\rho}{2}}+T^{\frac{3}{2}})\|F\|_{C([0,T];H^{-\rho}_x(\mathbb{R}))}.$$

Next, we want to establish estimates for the spatial traces in $H^1_t(0,T)$. Note that
\begin{equation}
\begin{aligned}W_t(x,t)&= \gamma_x(\partial_xS(t-t')[F(\cdot,t')])|_{t'=t} + \gamma_x\int_0^t \frac{d}{dt}S(t-t')[F_x(\cdot,t')]dt'    \\
	&=\int_{-\infty}^\infty ike^{ikx-k^2B(t)}\hat{F}(k,t)dk + \int_0^t b(t-t')\gamma_x\partial_x^3S(t-t')F(\cdot,t')dt'.
\end{aligned}
\end{equation}

So it follows 

\begin{equation}
\begin{aligned}
	\|W_t(x,\cdot)\|_{L^2_t(0,T)} &= \left(\int_0^T\left|\int_{-\infty}^\infty ike^{ikx-k^2B(t)}\hat{F}(k,t)dk\right|^2 dt\right)^{\frac12} \\
	& +  \left(\int_0^T\left|\int_0^tb(t-t')\gamma_x\partial_x^3S(t-t')F(\cdot,t')dt'\right|^2\right)^{\frac{1}{2}}\\
	&\le \int_0^T\left(\int_{-\infty}^\infty k^2e^{-2k^2B(t)}|\hat{F}(k,t)|^2dk\right)^{\frac12}dt \\
	& + \|b\|_{C[0,T]}\int_0^T\left(\int_0^t\left|\gamma_x\partial_x^3S(t-t')F(\cdot,t')\right|^2dt'\right)^{\frac12}dt\\
	& \lesssim (T^{\frac{3-\rho}{2}}+T^{\frac32})\|F\|_{C([0,T];H^{-\rho}_x(\mathbb{R}))} + (T^{\frac{1-\alpha}{2}}+T^{\frac{3}{2}})\|F\|_{C([0,T];H^{-\alpha}_x(\mathbb{R}))}, 
\end{aligned}
\end{equation}where $\alpha=1-\epsilon$ with $\epsilon>0$ arbitrarily small.  Combining $L^2$ and $H^1$ level estimates and interpolating, we establish the following lemma.
\begin{lemma}\label{nonhom_lem}
Let $-\frac{1}{2}\le s\le \frac{3}{2}$ and $F\in C([0,T];H_x^{s-\frac52+\epsilon}(\mathbb{R}))$, $\epsilon>0$ can be arbitrarily small, and $b\in C([0,T])$ with $b(t)\ge b_0>0$ for some $b_0>0$. Then, the solution of the non-homogeneous Cauchy problem \eqref{non_hom_cauchy} satisfies
the space estimate
\begin{equation}
	\sup_{t\in [0,T]}\|W(\cdot,t)\|_{H^s_x(\mathbb{R})}\lesssim \sqrt{T}\|F\|_{C([0,T];H^s_x(\mathbb{R}))}, \quad s\in \mathbb{R}.
\end{equation} and the time estimate
\begin{equation}
	\sup_{x\in \mathbb{R}} \|W(x,\cdot)\|_{H^{\frac{2s+1}{4}}_t(0,T)} \lesssim C_T \|F\|_{C([0,T];H_x^{s-\frac{5}{2}+\epsilon}(\mathbb{R}))},
\end{equation} where $C_T\rightarrow 0$ as $T\rightarrow 0^+$.
\end{lemma}

\begin{remark}
Again, noticing that we have $F_x$ at the right hand side of the equation instead of $F$, we can say the above lemma provides about $\frac{3}{2}$ units of improvement in temporal regularity of boundary traces compared to that of \cite[Theorem 3.2]{HMY19}.    This is due to parabolic smoothing. Such boost in regularity is for instance not expected for the Schrödinger equation.
\end{remark}

Note that a solution of \eqref{reduced_ibvp_q} can be obtained by finding a fixed point of the operator $\Upsilon(q)=q^0|_{(0,T)} + W|_{{[0,L]}},$ where $W$ depends on the input $q\in X_T := C([0,T];H_x^s(0,L))$ and solves the nonhomogeneous Cauchy problem with $F_x = c(t)\partial_x (Eq)$ (for a fixed bounded extension operator $E: H_x^s(0,L) \rightarrow H_x^s(\mathbb{R})$) and $q^0$ solves the reduced ibvp below that was studied in Case 1:

\begin{equation}\label{reduced_ibvp_q0}
\begin{aligned}
	q^0_t &= b(t)\partial_x^2q^0; \quad t\in (0,T');\\
	q^0(0) &= 0;\quad
	\gamma_0 q^0=(g-W(0,\cdot),-W(L,\cdot)),
\end{aligned}
\end{equation}

Considering the operator $\Upsilon$ on $X_T$, where $0\le s< \frac{3}{2}$, $s\neq \frac12$, in view of Lemma \ref{nonhom_lem} we have the estimate

\begin{equation}
\begin{aligned}
	\|\Upsilon(q)\|_{X_T}& \le \|q^0|_{(0,T)}\|_{X_T} + \|W|_{[0,L]}\|_{X_T}\\
	& \le \|g-W(0,\cdot)\|_{H_t^{\frac{2s+1}{4}}(0,T) } +\|W(L,\cdot)\|_{H_t^{\frac{2s+1}{4}}(0,T) } + \|W\|_{C([0,T];H_x^s(\mathbb{R}))}\\
	&\lesssim \|g\|_{H_t^{\frac{2s+1}{4}}(0,T) } + \|W(0,\cdot)\|_{H_t^{\frac{2s+1}{4}}(0,T)}+\|W(L,\cdot)\|_{H_t^{\frac{2s+1}{4}}(0,T) }+\sqrt{T}\|c(t)Eq\|_{C([0,T];H_x^s(\mathbb{R}))}\\
	&\lesssim \|g\|_{H_t^{\frac{2s+1}{4}}(0,T) } + C_T\|(c(t))Eq\|_{C([0,T];H_x^s(\mathbb{R}))}+\sqrt{T}\|c(t)Eq\|_{C([0,T];H_x^s(\mathbb{R}))}\\
	& \lesssim \|g\|_{H_t^{\frac{2s+1}{4}}(0,T) } + C_T\|q\|_{X_T}.
\end{aligned}
\end{equation}

In above estimate the constant $C_T$ now involves a factor of $\|c\|_{C[0,T]}$.  Now, let $R=2\|g\|_{H_t^{\frac{2s+1}{4}}(0,T) }$ where $T$ is small enough that $C_T<\frac{1}{2}$. Then, $\Upsilon$ maps the closed ball $B_R$ with radius $R$ in $X_T$ to itself. Furthermore, for any $q_1,q_2\in X_T$, we have 

\begin{equation}
\begin{aligned}
	\|\Upsilon(q_1)-\Upsilon(q_2)\|_{X_T}
	\lesssim  C_T\|q_1-q_2\|_{X_T},
\end{aligned}
\end{equation}namely $\Upsilon$ is a contraction on the same ball with small $T$. It follows that $\Upsilon$ has a unique fixed point in $B_R$. Furthermore, this fixed point is also unique in $X_T$. To see this assume existence of two fixed points $q, \tilde{q}\in X_T$. Then, $q-\tilde{q}$ solves the reduced ibvp with zero initial data and zero boundary data. Therefore, $q-\tilde{q}$ must be zero, namely $q=\tilde{q}$ in $X_T$.  Similarly, the solution of \eqref{reduced_ibvp_z} can also be obtained in view of the relation since it can be transformed into an ibvp of type \eqref{reduced_ibvp_q}.  Hence, combining these two foregoing results, we prove the following lemma for the case $c(t)\not\equiv 0.$

\begin{lemma}[Ibvp space estimate with time varying drift] \label{lwpforfull}
Let $0\le s< \frac32$, $s\neq \frac{1}{2}$, $c(t)\not\equiv 0$, for any $T>0$, $b,c \in C^{(\lceil s\rceil)}([0,T])$, $g_0,h_0\in H_t^{\frac{2s+1}{4}}(0,T)$, $u_0\in H_x^{s}(0,L)$ so that $U(0,\cdot), U(L,\cdot)\in H_t^{\frac{2s+1}{4}}(0,T)$.  If $s\in(\frac12,\frac32)$, assume that $(u_0,g_0)$ and $(u_0,h_0)$ satisfy the compatibility conditions $u_0(0)=g_0(0)$ and $u_0(L)=h_0(0)$, respectively. Then,  \eqref{reduced_ibvp} with $A(t)=b(t)\partial_{xx}+c(t)\partial_x$, has a unique solution on some interval $[0,T^*]$ for sufficiently small $T^*>0$ and it further satisfies the space estimate
\begin{equation}\label{space_est_y2}
	\sup_{t\in [0,T^*]}\|y(\cdot,t)\|_{H^s_x(0,L)}\lesssim \|g_0\|_{H_t^{\frac{2s+1}{4}}({0,T^*})} + \|h_0\|_{H_t^{\frac{2s+1}{4}}({0,T^*})} + \|u_0\|_{H_x^s(0,L)}
\end{equation} where the constant of the inequality only depends on the fixed parameters of the problem.

\end{lemma}

\subsection{Local and global well-posedness results for the full ibvp}
In this section, we state the main wellposedness results pertaining to  the cases of zero drift and nonzero drift, respectively. 

\begin{theorem}[Local and global solutions with zero drift]Let $T>0$, $s\ge 0$, $s\neq \frac{2n+1}{2}$ for $n\in \mathbb{N}$, $c(t)\equiv 0$, $b\in C^{(\lceil s\rceil)}([0,T])$, $b(t)\ge b_0$ for some $b_0>0$, $g_0,h_0\in H_t^{\frac{2s+1}{4}}(0,T)$, $u_0\in H_x^{s}(0,L)$.  Assume that $(u_0,g_0)$ and $(u_0,h_0)$ satisfy the necessary compatibility conditions at the space-time corner points $(x,t)=(0,0)$ and $(x,t)=(L,0)$, respectively. Then, \eqref{maineq} with $A(t)=b(t)\partial_{xx}$ has a unique solution $u\in C([0,T];H^s_x(0,L))$.  Moreover, the data to solution map $(u_0,g_0,h_0)\mapsto u$ is continuous from $H^s_x(0,L)\times \left[H_t^{\frac{2s+1}{4}}(0,T)\right]^2$ into $C([0,T];H^s_x(0,L))$.
\end{theorem}

\begin{proof}
Note that $u=U|_{(0,L)}+y|_{(0,T)},$ where $U$ solves the homogeneous Cauchy problem \eqref{hom_cauchy} with initial datum $U_0\in H^s_x(\mathbb{R})$ satisfying $U_0|_{(0,L)}=u_0$ and $\|U_0\|_{H^s_x(\mathbb{R})}\lesssim \|u_0\|_{H^s_x(0,L)}$ and $y$ solves the reduced ibvp \eqref{reduced_ibvp} with boundary data $g=g_0-U|_{x=0}$ and $h=h_0-U|_{x=L}$, extended in time with support in a larger time interval $(0,T')\supset (0,T)$. Therefore, the result follows by Lemma \ref{hom_lem}, Lemma \ref{hom_lem_time}, and Lemma \ref{zero_lem}.
\end{proof}

In above theorem, $T$ can be taken arbitrarily large assuming $b$ is globally defined and satisfies given conditions on each time interval $[0,T]$.  Namely, all solutions extend also globally in the case of zero drift.  In contrast, due to the fixed point argument, we first prove existence of solutions locally on a sufficiently small time interval in the case of nonzero drift. More precisely, we have the theorem below.

\begin{theorem}[Local solutions with nonzero drift]\label{local_nonzero}Let $0\le s< \frac{3}{2}$, $s\neq \frac12$, $c(t)\not\equiv 0$, for any $T>0$, $b,c\in C^{(\lceil s\rceil)}([0,T])$, $b(t)\ge b_0$ for some $b_0>0$, $g_0,h_0\in H_t^{\frac{2s+1}{4}}(0,T)$, $u_0\in H_x^{s}(0,L)$.  If $\frac32>s>\frac12$, assume that $(u_0,g_0)$ and $(u_0,h_0)$ satisfy the compatibility conditions $u_0(0)=g_0(0)$ and $u_0(L)=h_0(0)$, respectively. Then,  \eqref{maineq} with $A(t)=b(t)\partial_{xx}+c(t)\partial_x$ has a unique solution $u\in C([0,T^*];H^s_x(0,L))$ for some $T^*>0$ small enough.  Moreover, the data to solution map $(u_0,g_0,h_0)\mapsto u$ is continuous from $H^s_x(0,L)\times \left[H_t^{\frac{2s+1}{4}}(0,T^*)\right]^2$ into $C([0,T^*];H^s_x(0,L))$.
\end{theorem}

\begin{proof}
Again we set $u=U|_{(0,L)}+y|_{(0,T^*)}$, where $U$ solves the homogeneous Cauchy problem \eqref{hom_cauchy} with initial datum $U_0\in H^s_x(\mathbb{R})$ satisfying $U_0|_{(0,L)}=u_0$ and $\|U_0\|_{H^s_x(\mathbb{R})}\lesssim \|u_0\|_{H^s_x(0,L)}$ and $y$ is the solution of \eqref{reduced_ibvp}, this time with $c(t)\not\equiv 0$, on a time interval $[0,T^*]$ constructed in Lemma \ref{lwpforfull}.  Now, the result follows from Lemma \ref{hom_lem}, Lemma \ref{hom_lem_time}, Lemma \ref{nonhom_lem}, and Lemma \ref{lwpforfull}.
\end{proof}

\subsubsection*{Global solutions with nonzero drift}\label{Secglobal} Assume both $b$ and $c$ are globally defined, smooth, and uniformly bounded together with their derivatives with respect to $T$ on  all intervals $[0,T]$.
Theorem \eqref{local_nonzero} provides local well-posedness for the full linear convection-diffusion problem.  Local solutions can be proven to be global under suitable assumptions on data.  This can be done with the energy method at integer indexed spaces and then the result can be propagated towards fractional spaces.  Let us assume  without loss of generality $u_0\equiv 0$ and $h_0\equiv 0$ in \eqref{maineq} for simplicity since we can eliminate $u_0$ by subtracting from \eqref{maineq} the restriction to $(0,L)$ of the solution of the associated homogeneous Cauchy problem and then we can eliminate the resulting right hand boundary condition by splitting the problem into two linear problems which have zero boundary conditions at opposite endpoints and the second problem can be transformed into first one via a simple spatial change of variables.  Upon this reduction, i.e., assuming the only nonzero datum is the inhomogeneous left end Dirichlet one, multiplying the main equation by $u$ and integrating over $\Omega$, we get

\begin{equation}\label{firstmult}
\begin{aligned}
	\frac12\frac{d}{dt}\|u(\cdot,t)\|_{L^2_x(0,L)}^2 &= -b(t)\|u_x(\cdot,t)\|_{L_x^2(0,L)}^2 -b(t)u_x(0,t)g_0(t) + \frac{c(t)}{2}g_0^2(t)
\end{aligned}
\end{equation} The above identity involves the unknown Neumann trace at $x=0$. To gather information about this term, we  multiply the main equation by $u_{xx}$ and integrate over $\Omega$:

\begin{equation}\label{seciden}
\begin{aligned}
	\int_0^Lu_tu_{xx}dx=b(t)\int_0^Lu_{xx}^2dx+\frac{1}{2}c(t)u_x^2(L,t)-\frac{1}{2}c(t)u_x^2(0,t).
\end{aligned}
\end{equation}

Note that 

\begin{equation}
u_{x}^2(0,t) = -\frac{1}{L}\int_0^L[(L-x)u_x^2]_xdx\le \left(1+\frac{1}{L}\right)\|u_x(\cdot,t)\|_{L^2_x(0,L)}^2 + \|u_{xx}(\cdot,t)\|_{L^2_x(0,L)}^2.
\end{equation}

Similarly,

\begin{equation}
u_{x}^2(L,t) = \frac{1}{L}\int_0^L[xu_x^2]_xdx\le \left(1+\frac{1}{L}\right)\|u_x(\cdot,t)\|_{L^2_x(0,L)}^2 + \|u_{xx}(\cdot,t)\|_{L^2_x(0,L)}^2.
\end{equation}

The left hand side of \eqref{seciden} above can be written as

\begin{equation}\label{secondmult}
\int_0^Lu_tu_{xx}dx = -g_0'(t)u_x(0,t) -\frac{1}{2}\frac{d}{dt}\int_0^Lu_x^2(x,t)dx.
\end{equation}

Adding the identities \eqref{firstmult} and \eqref{secondmult}, we get

\begin{equation}
\begin{aligned}
	\frac12\frac{d}{dt}\|u(\cdot,t)\|_{H^1_x(0,L)}^2 & = -g_0'(t)u_x(0,t) -b(t)\|u_{xx}(\cdot,t)\|_{L^2_x(0,L)}^2 -\frac{1}{2}c(t)u_x^2(L,t)+\frac{1}{2}c(t)u_x^2(0,t)\\
	& -b(t)\|u_x(\cdot,t)\|_{L_x^2(0,L)}^2 -b(t)u_x(0,t)g_0(t) + \frac{c(t)}{2}g_0^2(t)\\
	&\le \frac{1}{\epsilon}(|g_0'(t)|^2+(b^2(t)+\frac{|c(t)|}{2})g_0^2(t))+(2\epsilon+|c(t)|-b(t))\|u_{xx}(\cdot,t)\|_{L^2_x(0,L)}^2\\
	&+(\frac{2(L+1)\epsilon}{L}+\frac{2(L+1)}{L}|c(t)|-b(t))\|u_{x}(\cdot,t)\|_{L^2_x(0,L)}^2.
\end{aligned}
\end{equation}  

It follows that if the sufficient condition $|c(t)|\le \frac{1}{2^+}\frac{L}{(L+1)}b(t)$ is satisfied for all $t\ge 0$ (i.e., if there is $\delta>0$, can be arbitrarily small, such that $|c(t)|\le\frac{1}{(2+\delta)}\frac{L}{(L+1)}b(t)$), then choosing $\epsilon>0$ small enough, we establish
the uniform bound on each time interval $[0,T]$ given by
\begin{equation}
\|u(\cdot,t)\|_{H^1_x(0,L)}\lesssim_{b,c,L}\|g_0\|_{H^1_t(0,T)}
\end{equation}
upon integrating the above identity and noticing that the coefficients of $\|u_{x}(\cdot,t)\|_{L^2_x(0,L)}^2$ and $\|u_{xx}(\cdot,t)\|_{L^2_x(0,L)}^2$ are negative so they can be dropped from the estimate. Namely, $H^1_x$-level local solutions can be extended globally in time assuming the modulus of the drift coefficient $c(t)$ is dominated by the diffusion coefficient $b(t)$ and $g_0$ is of sufficient regularity (here we had to assume $\frac14$ units more regularity for boundary input than what local theory suggested which was more optimal since it relied on local estimates whereas here we use the multiplier technique).  The following theorem follows from above arguments.
\begin{theorem} Let $b,c$ be globally defined and be as in Theorem \ref{local_nonzero} (with $s=1$) on each interval $[0,T]$, $g_0,h_0\in H_{t,\text{loc}}^1({0,\infty})$, $u_0\in H^1_x(0,L)$, $u_0(0)=g_0(0)$, $u_0(L)=h_0(L)$, and  assume in addition that $|c(t)|\le \frac{L}{2^+(L+1)}b(t)$ for $t\ge 0$, and $u\in C([0,T^*];H^1_x(0,L))$ be a local solution of \eqref{maineq} for some $T^*>0$ as in Theorem \ref{local_nonzero}.  Then, $u$ extends as a global solution to any time interval $[0,T]$.

\end{theorem}

\begin{remark}\label{rem-3.6}
In the zero-drift case, our analysis already shows that every local $H^s$ solution is in fact global. In the nonzero-drift case, however, the fixed-point approach yields local solutions only on small time intervals. The above theorem establishes global well-posedness under an additional assumption ensuring that the drift coefficient is suitably dominated by the diffusion coefficient. This restriction arises mainly from the multiplier technique we employ.

On the other hand, for the constant-coefficient model, the method used in the zero-drift case in Section \ref{Secredibvp} extends readily to the nonzero-drift problem. In particular, one obtains global well-posedness for the nonzero-drift model without the additional assumption in the above theorem. We therefore expect that a more direct approach based on the Fokas method that does not rely on the energy method, may weaken or even remove this assumption. This issue will be addressed in future work.

\end{remark}

\section{Other evolution PDEs}\label{Secevol}
In this section we follow the usual steps of the Fokas method, under the modified procedure described in section \ref{sec-sol-adv-dif} to derive the integral representation of the solution for IBVPs associated to some of the most well-known evolution PDEs of higher order.
\subsection{Linear KdV on the half-line}
Consider the PDE:
\begin{equation}\label{pde_kdv}
u_t + b(t) u_{xxx} + c(t) u_x = 0, \quad b(t) > b_0 > 0, \quad x>0, \ t>0,
\end{equation}
with the conditions
$
u(x,0)=u_0(x), \ \  x>0; \quad
u(0,t)= g_0(t), \ \ t>0.
$

\begin{figure}
\centering
\includegraphics[width=0.45\linewidth]{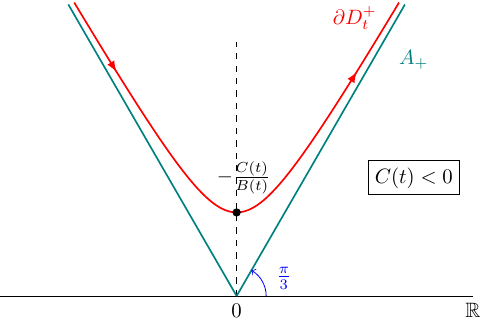} \qquad \includegraphics[width=0.45\linewidth]{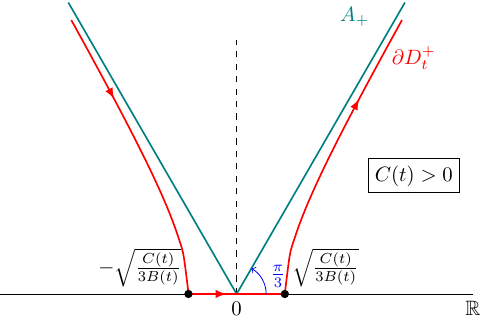}
\caption{The contours of integration $\partial D_t^+$.}
\label{fig:domain-kdv}
\end{figure}

We first introduce the half-line Fourier transform (HLFT) and its inverse, respectively as
\begin{equation}
\hat{\varphi}(k):=\int_0^\infty e^{-ikx}\varphi(x)dx,  \ \Im(k) \leq 0 ;\quad \varphi(x):=\frac{1}{2\pi}\int_{-\infty}^\infty e^{ikx}\hat{\varphi}(k)dk, \ x>0.
\end{equation}
Applying the HLFT in \eqref{pde_kdv} yields:
\begin{align*}
\hat{u}_t &  = -b(t) \int_0^\infty e^{-ikx} u_{xxx} \, dx - c(t) \int_0^\infty e^{-ikx} u_x \, dx \\ &= b(t) u_{xx}(0) + ikb(t) u_x(0) + b(t)(ik)^2  u(0)  + c(t) u(0) - b(t)(ik)^3 \hat{u} - c(t) ik \hat{u} \\
&= b(t) g_2(t) + ik b(t) g_1(t) - k^2 b(t) g_0(t) + c(t) g_0(t) + ik^3 b(t) \hat{u} - ik c(t) \hat{u} ,
\end{align*}
where $g_2$ and $g_1$ are unknown boundary conditions.
Thus,
\begin{equation}\label{kdv-help-GR}
\hat{u}_t + \left[-ik^3 b(t)  + ik c(t) \right] \hat{u}= b(t) g_2(t) + ik b(t) g_1(t) - k^2 b(t) g_0(t) + c(t) g_0(t). 
\end{equation}
Let \( B(t) = \int_0^t b(s) \, ds \) and \( C(t) = \int_0^t c(s) \, ds \), and define:
\begin{equation}\label{kdv-w-def}
w(k,t) = -ik^3 B(t) + ik C(t).
\end{equation}
By solving \eqref{kdv-help-GR}, we obtain 
\[ \hat{u}(k,t) e^{w(k,t)} = \hat{u}_0(k) + \int_0^t e^{w(k,s)} \left[ b(s)g_2(s) + ik b(s)g_1(s) - \left( k^2 b(s)-c(s) \right) g_0(s) \right] ds .\]
Defining the time-dependent transforms
$ g_j^f(k,t) = \int_0^t e^{w(k,s)} f(s) g_j(s) \, ds \, ,$
the global relation (GR) takes the form:
\begin{equation}
\hat{u}(k,t) e^{w(k,t)} = \hat{u}_0(k) + g_2^b(k,t) + ik g_1^b(k,t) - k^2 g_0^b(k,t)+ g_0^c(k,t), \qquad \Im(k) \leq 0 .
\end{equation}
Setting $G(k,t) =  g_2^b(k,t) + g_0^c(k,t)$, which is unknown, yields the following form of the GR
\begin{equation}\label{kdv-GR-1}
\hat{u}(k,t) e^{w(k,t)} = \hat{u}_0(k) +G(k,t) + ik g_1^b(k,t) - k^2 g_0^b(k,t), \qquad \Im(k) \leq 0 .
\end{equation}
By inverting the GR, and defining $D_t^+:=\{k\in\mathbb{C}\,|\,\text{Re}(\omega(k,t))<0, \ \Im k>0\}$ with $\partial D_t^+$ being the boundary of $D_t^+$ in the $k$-complex plane, we obtain the following integral representation of the solution:
\begin{align}\label{kdv-IR-1}
u(x,t)  = \int_\mathbb{R} e^{ikx-w(k,t)}\hat{u}_0(k)dk + \int_{\partial D_t^+}
e^{ikx-w(k,t)} \left[ G(k,t) + ik g_1^b(k,t) - k^2 g_0^b(k,t) \right] dk,
\end{align}
with $\partial D_t^+$ given in Figure \ref{fig:domain-kdv}.

The invariant mapping \( k \mapsto \nu(k) \) is derived from
\[ w(k,t)=w(\nu,t) \ \Leftrightarrow \   -ik^3 B(t) + ik C(t) = -i\nu^3 B(t) + i\nu C(t) \]
which simplifies to
\begin{equation}\label{kdv-nu}
\begin{aligned}
	(k - \nu)(k^2 B(t) + k\nu B(t) + \nu^2 B(t)- C(t)) = 0 \ \overset{\nu\ne k}{\Longrightarrow} \ \nu_{1,2}(k) = -\frac{k}{2} \pm \sqrt{-\frac{3}{4}k^2 + \frac{C(t)}{B(t)}} \\
	{\Longrightarrow} \ \nu_{1,2}(k) = -\frac{k}{2} \pm i \frac{\sqrt{3}}{2}\sqrt{k^2 - \frac{4}{3}\frac{C(t)}{B(t)}}.
\end{aligned}
\end{equation}

Evaluating the GR \eqref{kdv-GR-1} at \(\nu_1(k)\) and \(\nu_2(k)\) and {assuming wlog $C/B$ constant (again this assumption can be removed through a linear perturbation argument)}, we derive two additional global relations:
\begin{align*}
\hat{u}(\nu_1,t) e^{w(k,t)} &= \hat{u}_0(\nu_1) +  G(k,t) + i\nu_1  {g}_1^b(k,t) - \nu_1^2  {g}_0^b(k,t),  \\
\hat{u}(\nu_2,t) e^{w(k,t)} &= \hat{u}_0(\nu_2) +  G(k,t) + i\nu_2  {g}_1^b(k,t) - \nu_2^2  {g}_0^b(k,t),
\end{align*}
which are both valid for $ k \in D_t^+$, since $ k \in D_t^+ \Rightarrow \nu_j\in D_{j,t}^- \Rightarrow \Im \nu_j\le0, \ j=1,2$.

Solve the two GRs for \( {G}\) and \( {g}_1^b\), namely
\begin{align*}
& i {g}_1^b = ( \nu_1 + \nu_2)  {g}_0^b + \frac{1}{\nu_1-\nu_2} \left[ \hat{u}_0(\nu_2) - \hat{u}_0(\nu_1)- e^{w(k,t)}\left(\hat{u}(\nu_2,t) -\hat{u}(\nu_1,t)\right) \right],\\
&  {G} = -\nu_1 \nu_2   {g}_0^b + \frac{1}{\nu_1-\nu_2} \left[ \nu_1 \hat{u}_0(\nu_2) - \nu_2 \hat{u}_0(\nu_1)  - e^{w(k,t)}\left( \nu_1 \hat{u}(\nu_2,t) - \nu_2 \hat{u}(\nu_1,t) \right) \right]
\end{align*}
and using
$
\nu_1 + \nu_2 = -k \quad \text{and} \quad \nu_1 \nu_2 = k^2- \frac{C(t)}{B(t)},
$
we obtain
\begin{align*}
{G} + ik  {g}_1^b - k^2  {g}_0^b = & \left( -3k^2 + \frac{C(t)}{B(t)} \right)  {g}_0^b - \frac{\nu_1 - k}{\nu_1 - \nu_2} \hat{u}_0(\nu_2) - \frac{k - \nu_2}{\nu_1 - \nu_2} \hat{u}_0(\nu_1) \\
&+ e^{w(k,t)}\left[ \frac{\nu_1 - k}{\nu_1 - \nu_2} \hat{u}(\nu_2,t) + \frac{k - \nu_2}{\nu_1 - \nu_2} \hat{u}(\nu_1,t)  \right] .
\end{align*}

We substitute the above expression into the integral representation \eqref{kdv-IR-1} to find the solution:
\begin{align} \label{sol-kdv-1}
u(x,t) &= \frac{1}{2\pi} \int_{\mathbb{R}} e^{ikx - w(k,t)} \hat{u}_0(k) \, dk \notag \\
& - \frac{1}{2\pi} \int_{\partial D_t^+} e^{ikx - w(k,t)} \left[ \left( 3k^2 - \frac{C(t)}{B(t)} \right)  {g}_0^b + \frac{\nu_1 - k}{\nu_1 - \nu_2} \hat{u}_0(\nu_2) + \frac{k - \nu_2}{\nu_1 - \nu_2} \hat{u}_0(\nu_1) \right] dk,
\end{align}
where we have employed that $ \int_{\partial D_t^+} e^{ikx } \frac{\nu_j - k}{\nu_1 - \nu_2} \hat{u}(\nu_j,t) dk =0, \quad j=1,2. $

\subsubsection*{Special Case Reduction: $c(t) = 0$}
We obtain 
\[c(t) = 0 \ \Rightarrow C(t) = 0 \ \Rightarrow \quad w(k,t)=-ik^3 B(t) \ \Rightarrow \partial D_t^+=\partial D^+= \left\{ k = e^{i\frac{m\pi}{3}} r, \ r\ge 0, \ m=1,2  \right\} , \]
namely $\partial D^+$ is the boundary of $D^+=\left\{k: \ \arg k \in \left(\frac{\pi}{3},\frac{2\pi}{3} \right) \right\} $.
The terms in \eqref{sol-kdv-1} simplify as:
$\nu_j =e^{i\frac{2j\pi}{3}} k, \ j=1,2; \qquad
\frac{\nu_1 - k}{\nu_1 - \nu_2} =  e^{i\frac{4\pi}{3}} ; \qquad \frac{k - \nu_2}{\nu_1 - \nu_2} = e^{i\frac{2\pi}{3}} ;
$
hence,
\begin{align} \label{sol-kdv-c=0}
u(x,t) &= \frac{1}{2\pi} \int_{\mathbb{R}} e^{ikx - w(k,t)} \hat{u}_0(k) \, dk \notag \\
& - \frac{1}{2\pi} \int_{\partial D^+} e^{ikx - w(k,t)} \left[ 3k^2   {g}_0^b + e^{i\frac{4\pi}{3}}  \, \hat{u}_0\left(e^{i\frac{4\pi}{3}} k\right) + e^{i\frac{2\pi}{3}}  \, \hat{u}_0\left(e^{i\frac{2\pi}{3}} k\right) \right] dk.
\end{align}

\begin{remark}
In this subsection we derived the solution for $b(t)>b_0>0$. In the case of $b(t)<-b_0<0$, two boundary conditions are needed, e.g. $u(0,t)$ and $u_x(0,t)$; then one could modify the procedure of constant $b(t)=-1$, see \cite{DTV14}, and obtain the integral representation of the solution which involves the transforms of both boundary data.\\
We find interesting that if the sign of $b(t)$ is not fixed, then the above procedure proposes that the correct number of boundary conditions depends on the sign of $B(t)=\int_0^t b(s)ds$.
Indeed, the evolution feature of \eqref{pde_kdv} suggests that if $t>0$ belongs in an interval where 
\begin{itemize}
	\item $B(t)>0 \Rightarrow $ one boundary condition is required.
	\item $B(t)<0 \Rightarrow $ two boundary conditions are required.
\end{itemize}
\end{remark}

\subsection{Linear Schr\"odinger}

Consider the PDE on the finite interval
\begin{equation}\label{pde_schr}
iu_t + b(t)u_{xx} = 0, 
\end{equation}
with the conditions
\begin{align}
&u(x,0)=u_0(x), \ \  x\in(0,L) \\
&u(0,t)= g_0(t), \ \ u(L,t) = h_0(t), \ \  t>0.
\end{align}
Applying the FLFT and integrating by parts yields
\begin{align*}
i\hat{u}_t &= -b(t) \left[ e^{-ikx} u_x \big|_0^L + ik  e^{-ikx} u\big|_0^L +(ik)^2\hat{u} \right]\\
&= -b(t)\left[ h_1(t)e^{-ikL} - g_1(t) + ik \left(  h_0(t) e^{-ikL} - g_0(t) \right) - k^2\hat{u}\right]
\end{align*}
Solving this ODE for $\hat{u}$, we obtain the GR:
\begin{equation}\label{schr-GR}
e^{ik^2B(t)}\hat{u}(k,t) = \hat{u}_0(k) + ie^{-ikL}h_1^b - ig_1^b - ke^{-ikL}h_0^b + kg_0^b,\quad k\in\mathbb{C},
\end{equation}
where
$ f_j^b(k,t) = \int_0^t b(s) e^{i k^2 B(s)} f_j(s) \, ds, \ \  B(t) = \int_0^t b(s)\,ds$ .

The analogy with the case of constant $b(t)=b$ is obvious. Following the typical procedure as  in \cite{F02}, we obtain the following IR of solution:
\begin{align}\label{schr-IR}
u(x,t) &= \frac{1}{2\pi}\int_{\mathbb{R}} e^{ikx-ik^2B(t)}\hat{u}_0(k)dk \notag\\
& + \frac{1}{2\pi}\int_{\partial D^+} e^{ikx-ik^2B(t)}\left[ - ig_1^b + k g_0^b\right]dk \ + \frac{1}{2\pi}\int_{\partial D^-} e^{ik(x-L)-ik^2B(t)}\left[ -ih_1^b + kh_0^b \right]dk,
\end{align}
where $\partial D^\pm$ are the boundaries of the first and third complex $k$-quadrant, respectively.

Furthermore, the invariant mapping $k\to-k$ yields a second GR:
\begin{equation}\label{schr-GR-2}
e^{ik^2B(t)}\hat{u}(-k,t) = \hat{u}_0(-k) + ie^{ikL}h_1^b - ig_1^b + ke^{ikL}h_0^b - kg_0^b ,\quad k\in\mathbb{C}.
\end{equation}

Solving \eqref{schr-GR} and \eqref{schr-GR-2} for $h_1^b$ and $g_1^b$,  substituting the resulting expressions in \eqref{schr-IR} and simplifying we obtain the solution:
\begin{align}\label{schr-sol-1}
u(x,t) &= \frac{1}{2\pi}\int_{\mathbb{R}} e^{ikx-ik^2B(t)}\hat{u}_0(k)dk \notag \\
& + \frac{1}{2\pi}\int_{\partial \tilde{D}^+} \frac{e^{ikx-ik^2B(t)}}{e^{ikL}-e^{-ikL}}\left[2k h_0^b - 2k e^{-ikL}  g_0^b +e^{-ikL} \hat{u}_0(-k) -e^{ikL} \hat{u}_0(k)  \right]dk \\
& + \frac{1}{2\pi}\int_{\partial \tilde{D}^-} \frac{e^{ikx-ik^2B(t)}}{e^{ikL}-e^{-ikL}}\left[-2k h_0^b + 2k e^{-ikL}  g_0^b +e^{-ikL} \hat{u}_0(-k) -e^{-ikL} \hat{u}_0(k)  \right]dk, \notag
\end{align}
where $\partial \tilde{D}^\pm$ denote the slightly deformed $\partial D^\pm$ which avoid the roots of the denominator $e^{ikL}-e^{-ikL}$, namely the points $k_n=\frac{n \pi}{L},\ n\in\mathbb{Z}$ of the real line. 

The above expression can be also written in the following form
\begin{align*}\label{schr-sol-2}
u(x,t) &= \frac{1}{2\pi}\int_{\mathbb{R}} e^{ikx-ik^2B(t)}\hat{u}_0(k)dk+  \frac{1}{2\pi}\int_{\partial \tilde{D}^+} \frac{e^{-ik^2B(t)}}{\sin(kL)}2k \left[\sin(kx)h_0^b +\sin(k(L-x))g_0^b\right]dk  \notag \\
&  \ \ \frac{1}{2\pi}\int_{\partial \tilde{D}^+} \frac{e^{-ik^2B(t)}}{\sin(kL)}\left[\sin(kx) e^{ikL} \,\hat{u}_0(k) + \sin(k(L-x))\hat{u}_0(-k) \right]dk .
\end{align*}

\subsubsection*{Reduction in the case of half line}
Consider the PDE on the finite interval
\begin{equation}\label{pde_schr_hl}
iu_t + b(t)u_{xx} = 0, 
\end{equation}
with the conditions
$u(x,0)=u_0(x),  \ \ x>0;\quad u(0,t)= g_0(t), \ \ t>0.
$
One can imitate the procedure of \cite{F02}, under the above modifications to derive the solution of this IBVP. For matters of economy, using the uniform convergence of the associated integrals, we are allowed to take the limit $L\to\infty$, along with $h_0(t)=0$, in \eqref{schr-sol-1} to obtain the solution:
\begin{align}\label{schr-sol-3}
u(x,t) = \frac{1}{2\pi}\int_{\mathbb{R}} e^{ikx-ik^2B(t)}\hat{u}_0(k)dk 
+ \frac{1}{2\pi}\int_{\partial D^+} e^{ikx-ik^2B(t)} \left[2k   g_0^b - \hat{u}_0(-k)  \right]dk .
\end{align}

\subsection{Biharmonic Schr\"odinger}


Consider the PDE on the half line:
\begin{equation}\label{bih-hl}
i u_t +  b(t) u_{xxxx}  = 0, \quad x>0, \ t>0,
\end{equation}
with the conditions
$
u(x,0)=u_0(x), \ \  x>0; \quad
u(0,t)= g_0(t), \ \ u_{x}(0,t)= g_1(t), t>0.
$

Using the methodology of \cite{OY19} with the modifications of section \ref{sec-sol-adv-dif} we derive the solution
\begin{equation*}
\begin{aligned}
	u(x,t) = &\;
	\frac{1}{2\pi} \int_{-\infty}^{\infty} 
	e^{i k x + i k^4 B(t)} \, \hat u_0(k) \, dk 
	\\
	&- \frac{1}{2\pi} \int_{\partial D_1^+} 
	e^{i k x + i k^4 B(t)} \Big[ (1+i)\hat u_0(-i k) - i \hat u_0(-k) +
	2k^3(1-i)\, g_0^b
	+ 2k^2(1-i)\, g_1^b
	\Big] \, dk
	\\
	&- \frac{1}{2\pi} \int_{\partial D_2^+} 
	e^{i k x + i k^4 B(t)} 
	\Big[ (1-i)\hat u_0(i k) + i \hat u_0(-k) +
	2k^3(1+i)\, g_0^b
	- 2k^2(1+i)\, g_1^b
	\Big] \, dk,
\end{aligned}
\label{eq:modified_224}
\end{equation*}
where $\displaystyle g_j^b=\int_0^t b(s) g_j(s) e^{-ik^4 B(s)} ds$ and $\partial D_1^+$ and $\partial D_2^+$ are the boundaries of $D_1^+=\left\{k: \arg k\in \left(\frac{\pi}{4},\frac{\pi}{2}\right)  \right\}$ and $D_2^+=\left\{k: \arg k\in \left(\frac{3\pi}{4},\pi\right)  \right\}$, respectively.


\subsection{A fourth order PDE}

\subsubsection*{A half-line problem}
Consider the PDE on the half line:
\begin{equation}\label{pde_fourth-hl}
u_t + b(t) u_{xxxx}  = 0, \quad b(t) > b_0 > 0, \quad x>0, \ t>0,
\end{equation}
with the conditions
$
u(x,0)=0, x>0;\quad 
u(0,t)= g_0(t), u_{xx}(0,t)= g_2(t), t>0.
$
Applying the typical procedure, the global relation (GR) is
\begin{equation}\label{fourth-hl-gr}
e^{k^4 B(t)} \hat{u}(k,t) =   {g}^b_3(k,t) + ik   {g}^b_2(k,t) - k^2   {g}^b_1(k,t) - ik^3   {g}^b_0(k,t) , \qquad \Im k \le 0,
\end{equation}
where
\[   {g}^b_j(k,t) = \int_0^t b(s) e^{k^4 B(s)} g_j(s) \, ds, \quad B(s)=\int_0^t b(s) \, ds. \]
Inverting the GR, the integral representation of the solution reads
\begin{equation}\label{fourth-hl-ir}
u(x,t) = \frac{1}{2\pi} \int_{\partial D^+} e^{ikx - k^4 B(t)} \left[   {g}^b_3 + ik   {g}^b_2 - k^2   {g}^b_1 - ik^3   {g}^b_0 \right] dk,
\end{equation}
with $\partial D^+$ denoting the boundary $D^+$ consisting of the following two wedges
$$D^+=D_1^+ \cup D_2^+ =  \left\{ k: \ \arg k \in \left(\frac{\pi}{8}, \frac{3\pi}{8} \right)    \right\} \cup  \left\{ k:  \arg k \in  \left(\frac{5\pi}{8}, \frac{7\pi}{8} \right)    \right\}.$$
The invariant maps for $g_j^b(k,t)$ are $\nu_1(k)=ik, \ \nu_2(k)=-k, \nu_3(k)=-ik$, hence we  obtain two GR in $D_1^+$:
\begin{align*}
&g_3^b + (-ik) g_2^b + (-ik)^2 g_1^b + (-ik)^3 g_0^b= e^{k^4 B(t)} \hat{u}(-k,t) \\
&g_3^b + (-k) g_2^b + (-k)^2 g_1^b + (-k)^3 g_0^b= e^{k^4 B(t)} \hat{u}(ik,t) 
\end{align*}
and two GR in $D_2^+$:
\begin{align*}
&g_3^b + (-ik) g_2^b + (-ik)^2 g_1^b + (-ik)^3 g_0^b= e^{k^4 B(t)}  \hat{u}(-k,t) \\
&g_3^b + k g_2^b + k^2 g_1^b + k^3 g_0^b=e^{k^4 B(t)}  \hat{u}(-ik,t). 
\end{align*}
Solving for $ g_3^b$ and $ g_1^b $, each of the above systems, and substituting in \eqref{fourth-hl-ir} we obtain the solution
\begin{equation}\label{fourth-hl-sol}
u(x,t) = \frac{1}{2\pi} \int_{\partial D^+} e^{ikx - k^4 B(t)} 2 i k \left[    -{g}^b_2(k,t) + k^2    {g}^b_0(k,t) \right] dk.
\end{equation}

\subsubsection*{A finite-interval problem}
Consider the PDE on the finite interval:
\begin{equation}\label{pde_4th-hl}
u_t + b(t) u_{xxxx}  = 0, \quad b(t) > b_0 > 0, \quad x\in(0,L), \ t>0,
\end{equation}
with the conditions
\begin{align}
&u(x,0)=0, \ \  x\in(0,L) \\
&u(0,t) = g_0(t), \quad u_{xx}(0,t) = g_2(t), \ \ t>0,
&u(L,t) = h_0(t), \quad u_{xx}(L,t) = h_2(t), \ \ t>0.
\end{align}
Employing the usual methodology, one constructs the GR, which involve transforms of 4 known and 4 unknown boundary values; the associated IR of the solution is derived by the inversion of this GR. Then, producing 4 GRs valid in the whole complex $k$-plane, one can solve for the four unknown boundary values, and substitute in the IR, yielding the following solution :
\begin{align*}
u(x,t) = \frac{1}{2\pi} \int_{\partial D^+} e^{-k^4 B(t)} \frac{2ik}{\sin(kL)} \Big\{ \sin(kx) \left[ - {h}_2^b + k^2  {h}_0^b \right]  + \sin[k(L-x)] \left[ - {g}_2^b + k^2  {g}_0^b \right] \Big\} \, dk
\end{align*}
where $\partial D^+  =  \partial D_1^+  \cup  \partial D_2^+$ is defined as in the case of the half-line problem and
$ f_j^b(k,t) = \int_0^t b(s) e^{k^4 B(s)} f_j(s) \, ds . $

\section{Some remarks}\label{Secfinal}

\subsection{Connection to the half-line problem}

In section \ref{wpsection} we studied the well-posedness of the convection-diffusion equation on the finite interval. Of course, the special case $c(t)=0$ yields the analogue result for the heat equation which even improves the existing results for the case $b(t)=b>0$, constant. One shall expect that the well-posedness analysis for the heat equation on the half line, with time-variable coefficient $b(t)$, follows similar steps and yields equivalent regularity results. The derivation of estimates are indeed easier in the half-line case as it does not require analysis around $\lambda=0$. This suggests an alternative approach for proving regularity estimates  for the finite line problem because one can indeed reconstruct the solution of the finite line problem associated with the convection-diffusion equation by superposing solutions of two half-line problems on the finite line, one of which is posed on $(0,\infty)$ and a second one posed on $(-\infty,L)$. Therefore, regularity estimates for the finite line problem will directly follow from those estimates of the half-line problem.  Further details of this approach are given in \cite{KMO-pre}.

\subsection{Nonlinear problems}
The linear estimates derived in this paper can be applied to treat the local-well-posedness problem for the corresponding nonlinear problems, too. Consider for instance the nonlinear initial-boundary value problem
\begin{equation}\label{maineq_nonlin}
\begin{aligned}
	u_t = A(t)u + f(u); \quad t\in (0,T);\quad 
	u(0) = u_0;\quad
	\gamma_0 u=(g_0,h_0),
\end{aligned}
\end{equation}
where $f(u)=|u|^2u$ and $A(t)$ is defined as in \eqref{A_t}. Then, for fixed $z$ in a sought-after solution space, say $X_T=C([0,T];H_x^s(0,L))$, $s>\frac12$ replacing $f(u)$ by $f(z)$, one can first solve the following linear problem
\begin{equation}\label{maineq_lin_2}
\begin{aligned}
	u^z_t = A(t)u^z + f(z); \quad t\in (0,T);\quad 
	u^z(0) = u_0;\quad
	\gamma_0 u^z=(g_0,h_0),
\end{aligned}
\end{equation}
and then define the operator $\Upsilon: z\mapsto u^z$ on a closed ball, whose radius depend on relevant norms of initial-boundary data, in $X_T$. A fixed point of this operator for suitably chosen $T$, will then be solution to \eqref{maineq_nonlin}. The only necessary tools for this algorithm to yield a fixed point are the linear regularity estimates derived in the previous sections.  Technical details for construction of a fixed point through Banach fixed point theorem are almost identical to those given in \cite{HMY19} for the reaction-diffusion equation with constant coefficients, therefore omitted here.

\subsection{Robustness and new challenges}
The derivation of the integral representation for ibvps associated to evolution pdes with time-dependent coefficients in the current work demonstrates the unified character of the Fokas method. We consider the amendments to the typical approach (with constant coefficients) direct, with the most important one to be the introduction of a time-dependent contour of (spectral) integration. The fact that these contours are asymptotically independent from $t$ contributes to the numerical applicability of the method. These arguments are valid for all the evolution equations considered here.

As mentioned in Remark \ref{rem-3.6}, time dependent coefficients in the framework of multiple spatial derivatives bring some challenges for direct estimation of the associated regularity. This was partially overcome by the employment of fixed point and multiplier arguments. New insights on the direct usage of the explicit integral representation formula are needed for more general regularity results.

\enlargethispage{20pt}
\section*{Acknowledgements}{KK acknowledges support by the Sectoral Development Program (SDP 5223471) of the Ministry of Education, Religious Affairs and Sports, through the National Development Program (NDP) 2021-25, grant no 200/1029.}


\bibliographystyle{plain}
\bibliography{references}

\end{document}